\documentclass[11pt]{amsart}

\usepackage{amssymb}
\usepackage{amscd}

\newcommand{\nc}{\newcommand}

\newtheorem{theorem}{Theorem}[section]
\newtheorem{prop}[theorem]{Proposition}
\newtheorem{importnota}[theorem]{Important Notation}
\newtheorem{prblm}[theorem]{Problem}
\newtheorem{notation}[theorem]{Notation}
\newtheorem{caution}[theorem]{Caution}
\newtheorem{remark}[theorem]{Remark}
\newtheorem{lemma}[theorem]{Lemma}
\newtheorem{construction}[theorem]{Construction}
\newtheorem{corollary}[theorem]{Corollary}
\newtheorem{example}[theorem]{Example}
\newtheorem{conclusion}[theorem]{Conclusion}
\newtheorem{triviality}[theorem]{Triviality}
\newtheorem{proto}[theorem]{Prototype Quasifibration}
\newtheorem{cauex}[theorem]{Cautionary Example}
\newtheorem{propositiondef}[theorem]{Proposition-Definition}
\newtheorem{subth}{Nuisance}[theorem]
\newtheorem{ssubth}{ }[subth]
\newtheorem{conjecture}[theorem]{Conjecture}
\newtheorem{sidest}[theorem]{Side Story}
\newtheorem{miniexample}[theorem]{Example}
\theoremstyle{definition}
\newtheorem{defin}[theorem]{Definition}

\nc\tri[1]{\begin{triviality}}
\nc\side[1]{\begin{sidest}}
\nc\conj[1]{\begin{conjecture}}
\nc\prodef[1]{\begin{propositiondef}}
\nc\prt[1]{\begin{proto}}
\nc\lem[1]{\begin{lemma}}
\nc\sblm[1]{\begin{sublemma}}
\nc\pro[1]{\begin{prop}}
\nc\thm[1]{\begin{theorem}}
\nc\cor[1]{\begin{corollary}}
\nc\dfn[1]{\begin{defin}}
\nc\sthm[1]{\begin{subth}}
\nc\exm[1]{\begin{example}}
\nc\miniexm[1]{\begin{miniexample}}
\nc\plm[1]{\begin{prblm}}
\nc\rmk[1]{\begin{remark}}
\nc\subrmk[1]{\begin{subremark}}
\nc\ntn[1]{\begin{notation}}
\nc\cau[1]{\begin{caution}}
\nc\imn[1]{\begin{importnota}}
\nc\cax[1]{\begin{cauex}}
\nc\con[1]{\begin{construction}}
\nc\ssthm[1]{\begin{ssubth}}
\nc\cnc[1]{\begin{conclusion}}
\nc\elem{\end{lemma}}
\nc\esblm{\end{sublemma}}
\nc\eside{\end{sidest}}
\nc\econj{\end{conjecture}}
\nc\eprodef{\end{propositiondef}}
\nc\eprt{\end{proto}}
\nc\ethm{\end{theorem}}
\nc\ecor{\end{corollary}}
\nc\edfn{\end{defin}}
\nc\esthm{\end{subth}}
\nc\epro{\end{prop}}
\nc\etri{\end{triviality}}
\nc\eexm{\end{example}}
\nc\eminiexm{\end{miniexample}}
\nc\ermk{\end{remark}}
\nc\subermk{\end{subremark}}
\nc\eplm{\end{prblm}}
\nc\ecau{\end{caution}}
\nc\ecax{\end{cauex}}
\nc\eimn{\end{importnota}}
\nc\entn{\end{notation}}
\nc\econ{\end{construction}}
\nc\ecnc{\end{conclusion}}
\nc\essthm{\end{ssubth}}

\newcommand{\C}{\mathbb{C}}
\newcommand{\R}{\mathbb{R}}
\newcommand{\Q}{\mathbb{Q}}

\newcommand{\Z}{\mathbb{Z}}

\newcommand{\N}{\mathbb{N}}
\newcommand{\F}{\mathbb{F}}
\newcommand{\A}{\mathbb{A}}

\newcommand{\p}{\mathfrak{p}}
\newcommand{\w}{\varpi}
\newcommand{\di}{{\rm diag}}

\renewcommand{\o}{\mathfrak{o}}
\newcommand{\ds}{\displaystyle}

\newcommand{\lra}{\longrightarrow}
\newcommand{\one}{\mathbf{1}_{E^1}}

\begin{document}

\title[$p$-stabilized vectors]{An explicit computation of $p$-stabilized vectors}

\author{\sc Michitaka MIYAUCHI}
\address{Faculty of Liberal Arts and Sciences\\
Osaka Prefecture University\\
1-1 Gakuen-cho, Nakaku, Sakai, Osaka 599-8531, JAPAN}
\email{michitaka.miyauchi@gmail.com}

\author{\sc  Takuya YAMAUCHI}
\address{
Department of mathematics, Faculty of Education\\
Kagoshima University\\
Korimoto 1-20-6 Kagoshima 890-0065, JAPAN  and\\ 
Department of mathematics \\
 University of Toronto \\
Toronto, Ontario M5S 2E4, CANADA}
\email{yamauchi@edu.kagoshima-u.ac.jp or tyama@math.toronto.edu}


\maketitle


\begin{abstract}
In this paper, 
we give a method to compute $p$-stabilized vectors in the space of 
parahori-fixed vectors for connected reductive groups over $p$-adic fields. 
An application to the global setting is also discussed. 
In particular, we give an explicit $p$-stabilized form of a Saito-Kurokawa lift. 
\end{abstract}

\maketitle
\pagestyle{myheadings}
\markboth{}{}

\section{Introduction}
Let $F$ be a non-archimedean local field, 
$\o$ the ring of integers of $F$, 
$\p$ the maximal ideal of $\o$, 
$\w$ a uniformizer of $F$, and 
$\F=\o/\p$ the residue field of $F$. 
We normalize the valuation $|\cdot |$ of $F$ 
so that $|\w|=q^{-1}$, where $q$ is the cardinality of $\F$.  
Let $G$ be a connected reductive group defined over $F$, 
$B$ its standard Borel subgroup 
and $K=G(\o)$ a maximal compact subgroup of 
$G(F)$ whenever it is defined. 
Let $P$ be a parabolic subgroup of $G$ which contains $B$ and 
$K_P=\{g\in K|\ g\ {\rm mod}\ \p\in P(\F)\}$ 
the parahoric subgroup of $K$ associated to $P$. 

Let $\pi$ be an irreducible smooth representation of $G(F)$ 
such that the space $\pi^{K_P}$ of $K_P$-fixed vectors 
in $\pi$ is non-trivial.
Then the Hecke algebra $\mathcal{H}_{K_P}$ 
of $G(F)$
associated to  $K_P$ acts on $\pi^{K_P}$. 
In this paper, 
we first give a method to compute the eigenvalues for
the special elements of 
$\mathcal{H}_{K_P}$ on $\pi^{K_P}$, which are called ``$U_p$-operators". We next give an explicit construction of 
 simultaneous eigenvectors  for these $U_p$-operators,
 which are called ``$p$-stabilized vectors".   

An idea to compute eigenvalues of $U_p$-operators
is to consider the Jacquet module of $\pi$ associated to $P$. When $P$ is the standard Borel subgroup of $G(F)$, 
this has been a well-known method for experts. 
We extend  a result  of Casselman \cite{cas-note} which is proved only for the standard Borel subgroup 
to any parabolic subgroup
$P$ and any compact open subgroup of $G(F)$ 
which has an Iwahori factorization relative to $P$ in the sense of \cite{cas-note}. 
Then it can be reduced a calculation of 
the eigenvalues for $U_p$-operators 
to the same problem for the actions of 
the elements corresponding to those operators
on the Jacquet modules (Proposition~\ref{plus}).

Our results might give a potential tool to study an arithmetic investigation of automorphic forms in 
Iwasawa theory, Hida theory, or  deformation theory of Galois representations 
\cite{Mazur},\cite{til},\cite{skinner&urban}, though we do not discuss about this in this paper. 
In particular, we will know that what kind of $p$-stabilized forms can be embedded into a Hida family with respect a specific 
parabolic subgroup of $G$.  

This paper is organized as follows. 
In Section 2, we study the action
of Hecke algebras 
on the space of parahori-fixed vectors
by using Jacquet modules.  
In Section 3, we introduce the notion of  $U_p$-operators and $p$-stabilized vectors. 
In Section 4, we outline a method to construct  $p$-stabilized vectors when  $\pi$ has 
a non-zero $K$-fixed vector. Without this assumption  on $\pi$, it seems to be difficult to check the non-triviality of the vector 
which we will construct. 
In Section 5, 6, and 7, we make up a list of all $U_p$-eigenvalues and $p$-stabilized vectors 
for $GL_2, U(2,1)$, and $GSp_4$. 
A relation to the global setting is discussed in Section 8 and then 
the global $p$-stabilized forms are given in the final section in cases of $GL_2$ and $GSp_4$. 
In particular, we will give them for Saito-Kurokawa lifts where the existence has been already discussed in 
Proposition 4.2.2, p.688 of \cite{skinner&urban} (see (6) and (7)).    

The authors would like to thank Professor Tomonori Moriyama
for helpful conversations. The second author
is partially supported by JSPS Grant-in-Aid for Scientific Research No.23740027 and JSPS Postdoctoral 
Fellowships for Research Abroad No.378.

\section{The action of Hecke operators via Jacquet modules}
We keep the notation in Section 1. 
In this section, we study parahori-fixed vectors of 
smooth representations of $G(F)$. 
We fix a maximal torus $T$ of $G$ and a minimal parabolic subgroup $B$ of $G$ which contains $T$. 
Then we have the Levi decomposition $B = TU$,
where $U$ is the unipotent radical of $B$. 
Let $P$ be a parabolic subgroup of $G$ containing $B$ with  Levi decomposition $P=MN$.
We denote by
$\overline{N}$ the unipotent radical of the parabolic subgroup opposite to $P$. 

Henceforth, for any algebraic group $H$, we sometimes denote by $H$ the group of $F$-valued points of $H$ 
for the sake of simplicity. This should cause no confusion in the remainder of the paper.
   
For any smooth representation $(\pi,V)$ of $G(F)$,
we define
its Jacquet module $(\pi_N,V_N)$ as follows (cf. Section 3 of \cite{cas-note}):
Set $V(N) =  \langle \pi(n) v-v\, |\, n\in N,\ v\in V
\rangle$
and $V_N = V/V(N)$.
We define a representation $(\pi_N, V_N)$ of $M$ by 
$$\pi_N(m)r_N(v)=\delta^{-\frac{1}{2}}_P(m)
r_N(\pi(m)v),\ m\in M,\ v\in V,$$
where $r_N$ is the natural projection 
from $V$ to $V_N$
and
$\delta_P$ is the modulus character of $P(F)$. 
 
Let $I:=\{g\in K\, |\, g\bmod \p\in B(\o/\p)\}$ be the
standard Iwahori subgroup of $G(F)$. 
If a smooth representation $(\pi, V)$ of $G(F)$ is admissible, then by p.7, Theorem of \cite{garrett}, the canonical projection $r_U: V\rightarrow  V_U$ 
induces a $\C$-linear isomorphism
\begin{eqnarray}
V^I\stackrel{\sim}{\lra}(V_U)^{I\cap T(F)}.
\end{eqnarray}
The following theorem generalizes this isomorphism 
to any parabolic subgroups.
\thm{}\label{key-thm} 
Let $(\pi, V)$ be an admissible representation of $G(F)$.
Suppose that a compact subgroup $J$ of $G(F)$ 
contains $I$
and it has an Iwahori factorization
with respect to $P$ 
in the sense of \cite{cas-note} (see before Proposition 1.4.4 in loc.cit.).
Then the canonical projection $r_N:V\rightarrow V_N$ induces 
a linear isomorphism $$V^J\stackrel{\sim}{\lra} (V_N)^{J\cap M}.$$  
\ethm
\begin{proof}
By Theorem 3.3.3 of \cite{cas-note}, the map $r_N:V^J\lra (V_N)^{J\cap M}$  is surjective. 
We now prove the injectivity of this map. 
By (1), we get
$V^I\cap V(U)=\{0\}$. 
Since $I\subset J$ and $U\supset N$, we have 
$V^I\supset V^J$ and $V(U)\supset V(N)$. This gives us that $V^J\cap V(N)\subset V^I\cap V(U)=\{0\}$. 
This implies that $r_N:V^J\lra (V_N)^{J\cap M}$ is 
injective.
\end{proof}

As in Theorem~\ref{key-thm},
let $J$ be a compact subgroup of $G(F)$ which contains $I$. Assume that $J$ has
an Iwahori factorization with respect to $P$.

\dfn. (\cite{bk} Definition 6.5)\label{plus-minus}
We say that an element $m\in M$ is {\it positive
relative to} $(P,J)$ 
if the following conditions are fulfilled: 
$$m(J\cap N) m^{-1}\subset J\cap N,\ m^{-1}(J\cap \overline{N}) m
\subset J\cap \overline{N}.$$ 
We denote by $M^+$ the set of all positive elements in $M$.
We say that an element $m$ in $M$ 
is  {\it negative relative to} $(P,J)$ 
if $m^{-1}$ is positive.
We write $M^-$ for the set of all negative elements in $M$.
\edfn

Given a compact open subgroup $J$ of $G(F)$,
we define the Hecke algebra 
$\mathcal{H}_J:=\mathcal{H}[G(F)//J]$
of $G(F)$ associated to $J$ to be the space of all 
compactly supported functions $f:G(F)\lra \C$ 
which satisfy
$f(j_1gj_2)=f(g)$, for $j_1,j_2\in J$ and $g\in G(F)$. 
Then $\mathcal{H}_J$ becomes an algebra 
under the convolution
with respect to the Haar measure on $G(F)$
normalized so that the volume of $J$ is one.
For any $g \in G(F)$,
we denote by $f_g=[JgJ]\in \mathcal{H}_J$ 
the characteristic function of $Jg J$. 
Since the algebra $\mathcal{H}_J$ is generated by $f_g,\ g\in G(F)$, if we define 
$\mathcal{H}_{J,\Q}:=\Q[f_g\,|\, g\in G(F)]$,
 then we have $\mathcal{H}_J=\mathcal{H}_{J,\Q}\otimes_\Q\C$. 
For any $\Q$-algebra $A$,
we put $\mathcal{H}_{J,A}:=\mathcal{H}_{J,\Q}\otimes_\Q A$.

If $(\pi, V)$ is a smooth representation of $G(F)$,
then the Hecke algebra $\mathcal{H}_J$ acts on 
$V^J$ (cf. \cite{cas-note}). 
We denote by $Z_M$ the center of $M$. 
We consider 
Hecke operators associated to positive elements
in $Z_M$. 
\pro{}\label{plus} Let $J$ be as in Definition \ref{plus-minus} and $(\pi, V)$ an
admissible representation of $G(F)$. 
Then for any $\zeta\in Z_M\cap M^+$,
we have
$$
r_N(\pi(f_\zeta)v)
= \delta^{-\frac{1}{2}}_P(\zeta)\pi_N(\zeta)r_N(v),\
v \in V^J,
$$
where 
$r_N:V^J\stackrel{\sim}{\lra}(V_N)^{J\cap M}$
is the isomorphism given in Theorem~\ref{key-thm}. 
\epro
\begin{proof}
For $\zeta\in Z_M\cap M^+$ and $v\in V^J$,
we have 
$$\pi(f_\zeta)v=\int_{J\zeta J}\pi(g)vdg=\sum_{k\in J/J\cap \zeta J\zeta^{-1}}\pi(k\zeta)v.$$
Since we assume that $J$ has an Iwahori factorization,
we get
$J=(J\cap \overline{N})(J\cap M)(J\cap N)$.
Because
$\zeta$ is positive and it belongs to $Z_M$, 
we have that
$J\cap \zeta J\zeta^{-1}
= (J\cap \overline{N})(J\cap M)\zeta (J\cap N)\zeta^{-1}$,
so that
$J/J\cap \zeta J\zeta^{-1}
= (J\cap N)/\zeta (J\cap N)\zeta^{-1}$.
Therefore we obtain
$$\pi(f_\zeta)v=\sum_{k\in (J\cap N)/\zeta (J\cap N)\zeta^{-1}}\pi(k\zeta)v.$$
Since $r_N(\pi(k\zeta)v) = r_N(\pi(\zeta)v)
= \delta^{\frac{1}{2}}_P(\zeta)\pi_N(\zeta)r_N(v)$,
$k \in J\cap N$
and 
$[J\cap N:\zeta (J\cap N)\zeta^{-1}]
= \delta^{-1}_P(\zeta)$,
we have 
$$r_N(\pi(f_\zeta)v)
=\delta^{-\frac{1}{2}}_P(\zeta)\pi_N(\zeta)r_N(v),$$
as required.
\end{proof}

\pro{}\label{commutativity} 
Suppose that a compact open subgroup $J$
of $G(F)$ has an Iwahori factorization relative
to $P$.
Then
$f_{\zeta_1}$ and $f_{\zeta_2}$ are
commutative,
for any $\zeta_1$, $\zeta_2 \in Z_M\cap M^+$.
\epro
\begin{proof}
We shall claim that 
$f_{\zeta_1}*f_{\zeta_2} = f_{\zeta_1 \zeta_2}$.
Then we obtain
\[
f_{\zeta_1}*f_{\zeta_2} = f_{\zeta_1 \zeta_2}
= f_{\zeta_2 \zeta_1}=
f_{\zeta_2}*f_{\zeta_1}
\]
because $\zeta _1 \zeta_2 =\zeta_2\zeta_1$.
It follows from \cite{bk} (6.6) that
\[
[J\zeta_1 J: J]
= [J\cap \overline{N}: \zeta_1^{-1}(J\cap \overline{N})\zeta_1]
[J\cap M: \zeta_1^{-1}(J\cap M)\zeta_1\cap (J\cap M)].
\]
Hence we get 
\[
[J\zeta_1 J: J]
= [J\cap \overline{N}: \zeta_1^{-1}(J\cap \overline{N})\zeta_1]
\]
since $\zeta_1$ lies in the center of $M$.
Similarly,
we obtain 
$[J\zeta_2 J: J]
= [J\cap \overline{N}: \zeta_2^{-1}(J\cap \overline{N})\zeta_2]$
and
$[J\zeta_1\zeta_2 J: J]
= [J\cap \overline{N}: (\zeta_1\zeta_2)^{-1}(J\cap \overline{N})\zeta_1\zeta_2]$.
Since $\zeta_1$ and $\zeta_2$ are both positive,
we have
\[
J\cap \overline{N} \supset 
\zeta_2^{-1}(J\cap \overline{N})\zeta_2
\supset 
(\zeta_1\zeta_2)^{-1}(J\cap \overline{N})\zeta_1\zeta_2.
\]
So we obtain
\[
[J\zeta_1\zeta_2 J: J]
= [J\zeta_1J: J] [J\zeta_2 J: J],
\]
and hence
$f_{\zeta_1}*f_{\zeta_2} = f_{\zeta_1 \zeta_2}$ by Proposition 2.2 in Chapter 3 of \cite{Harish}.
\end{proof}

\section{$p$-stabilized vectors}
For simplicity, we assume that the dimension of the center of $G$ is at most one.  
Let $\pi$ be an irreducible smooth representation of $G(F)$. 
Then $\pi$ is admissible by \cite{jacquet}. 
In this section, we  give a notion of $p$-stabilized vectors (or of $p$-stabilization ) for 
parahori-fixed vectors in $\pi$. 

Let $\Delta$ be the set of all simple roots of $(G,T)$ which is a subset of the character group 
$X^\ast(T):={\rm Hom}_{{\rm alg}}(T,GL_1)$. 
Let $P$ be a parabolic subgroup of $G$ containing $B$, 
$P=MN$ its Levi decomposition, 
and $K_P$ the parahoric subgroup which corresponds to $P$. 
Let $\Delta_P$ be the subset of $\Delta$ corresponding to $P$. 
We define $T^-_P$ to be the semi-group consisting of the elements $t$ in
$T(F)/T(\o)$ such that 
\begin{eqnarray}
|\alpha(t)|\le 1 \mbox{ for all $\alpha\in \Delta$ and $t(K_P\cap N)t^{-1}\subset K_P\cap N$}.
\end{eqnarray}
We can choose a complete system of representatives for $T^-_P$ as elements in $Z_M\cap M^+$. 
Put $m_P=\sharp(\Delta\setminus\Delta_P)$. Note that $\Delta_B=\emptyset$.  
For each $\alpha \in\Delta$, there exists $t_\alpha\in T^-_B$ such that 
${\rm ord}_\w\alpha(t_\alpha)=1$ and ${\rm ord}_\w\alpha(t_\beta)=0$ for all $\beta\in \Delta\setminus\{\alpha\}$.
Put 
$$t_{m_P+1}:=\left\{
\begin{array}{cc}
\w^{-1}{\rm Id}, & \mbox{if $Z_{M}\supset F^\times$}\\
{\rm Id}, & \mbox{otherwise},
\end{array}
\right.
$$
where Id is the identity element of $G(F)$. 
We write $\Delta\setminus\Delta_P=\{\alpha_1,\ldots,\alpha_{m_P}\}$ and $t_i=t_{\alpha_i}$ for $i\in\{1,\ldots,m_P\}$. 
Then $T^-_P$ is generated by $t_1,\ldots,t_{m_P},t_{m_P+1}$ as a semi-group.
For any $\Q$-algebra $A$ which is contained in $\C$, 
we consider the subalgebra $$\mathcal{U}_{P,A}:=A[[K_PtK_P]\ |\ t\in 
T^-_P ]$$ of the Hecke algebra $\mathcal{H}_{K_P,A}$ over $A$. 
\lem{}\label{commutative} 
Put $U^P_{\w,i}:=[K_Pt_iK_P]\in \mathcal{H}_{K_P}$,
for $i\in\{1,\ldots, m_P+1\}$. 
Then the ring $\mathcal{U}_{P,A}$ is 
a commutative $A$-algebra generated by  
$U^P_{\w,1},\ldots,U^P_{\w,m_P+1}$. 
\elem 
\begin{proof}
Recall that we take a complete system of representatives for $T_P^-$ as elements in $Z_M\cap M^+$. 
Then the commutativity follows from Proposition \ref{commutativity}. The later claim follows from the fact that  
$T^-_P$ is generated by $t_1,\ldots,t_{m_P},t_{m_P+1}$.
\end{proof}
\dfn{} Let $(\pi, V)$ be an irreducible smooth 
representation of $G(F)$ such that $V^{K_P}\neq \{0\}$. 
We say that a non-zero vector $v$ in $V^{K_P}$ is {\it a $p$-stabilized vector with respect to} $\widehat{P}$ if 
it is a  simultaneous eigenvector for all 
$U_{\w,1}^P,\ldots,U_{\w,m_P}^P$. Here $\widehat{P}$ is the Langlands dual of $P$ (cf. \cite{Borel}).  
\edfn
\rmk{}
The condition (2) on $T_P^-$ is crucial to get the commutativity of $\mathcal{U}_{P,A}$. 
In general, this property does not hold for $\mathcal{H}_{K_P,A}$. 
\ermk

\section{Construction of $p$-stabilized vectors}
Let $(\pi, V)$ be an irreducible smooth representation of 
$G(F)$ which has a non-zero $K$-fixed vector. 
In this section, we give a method
to produce $p$-stabilized vectors for $\pi$.
Let $P=MN$ be a parabolic subgroup of $G(F)$
containing $B$. 
Then by Theorem~\ref{key-thm},
the Jacquet functor $r_N$ induces an isomorphism
$r_N: V^{K_P} \simeq (V_N)^{K_P\cap M}$.
We set $W = (V_N)^{K_P\cap M}$.
Let $H$ denote the subgroup of $Z_M$
generated by $t_1, \ldots, t_{m_P} \in T_P^-$. As we have seen before, 
we may assume $H\subset Z_M\cap M^+$. 
For a quasi-character $\chi$ of $H$
and $n \in \N$,
we define
\[
W_{\chi, n}
=\{w \in W\, |\, (\pi_N(t)-\chi(t))^n w = 0 \ \mbox{for any}\  
t \in H\}
\]
and put $W_{\chi,\infty}=\bigcup_{n \in \N}W_{\chi,n}$. 
Similarly, we define $(V_N)_{\chi, \infty}$ for $V_N$.
Let $\mathcal{S}$ denote the set of 
quasi-characters $\chi$ of $H$
such that $W_{\chi, \infty} \neq \{0\}$.
Since $W$ is a finite-dimensional $H$-module,
we have
\[
W = \bigoplus_{\chi \in \mathcal{S}} W_{\chi, \infty}.
\]
For an element $w$ in $W$,
$w$ is a simultaneous eigenvector for $t_1, \ldots, t_{m_P}$
if and only if 
$w$ lies in $W_{\chi,1}$, for some $\chi \in \mathcal{S}$.

Let $\phi_K$ be a non-zero $K$-fixed vector in $V$.
By the Iwasawa decomposition $G = PK$,
the element $v = r_N(\phi_K)$ generates
$V_N$ as an $M$-module,
so does $W$.
Since $H$ is contained in  the center of $M$,
we have
\[
V_N = \bigoplus_{\chi \in \mathcal{S}} (V_N)_{\chi, \infty}
\]
as an $M$-module.
We claim that the $W_{\chi, \infty}$-component of $v$
is not zero, for any $\chi \in \mathcal{S}$.
If the $W_{\chi, \infty}$-component of $v$
is zero,
then 
$v$ lies in the proper $M$-submodule
$ \bigoplus_{\chi' \neq \chi} (V_N)_{\chi', \infty}$
of $V_N$.
This contradicts the fact that 
 $v$ generates
$V_N$ as an $M$-module.
So the claim follows.

We fix a character $\chi$ of $H$ in $\mathcal{S}$.
For any $\chi' \in \mathcal{S}$
which is different from $\chi$,
there exists an integer $1 \leq i(\chi') \leq m_P$
such that $\chi(t_{i(\chi')}) \neq \chi'(t_{i(\chi')})$.
Put $n(\chi') = \dim W_{\chi', \infty}$.
Then 
\[
\displaystyle 
v' = 
\prod_{\chi' \neq \chi} (\pi_N(t_{i(\chi')})-\chi'(t_{i(\chi')}))^{n(\chi')} v
\]
is a non-zero vector in $W_{\chi, \infty}$.
Therefore,
there exist non-negative integers 
$n(\chi, i)$ for $1 \leq i \leq m_P$
such that
$v'' := \prod_{1\leq i\leq m_P}(\pi_N(t_{i})-\chi(t_{i}))^{n(\chi, i)}v' $
is a non-zero vector
in $W_{\chi,1}$.
By Proposition~\ref{plus},
\[
\displaystyle 
\phi =
\prod_{1\leq i \leq m_P}
(\delta_P^{\frac{1}{2}}(t_i)\pi( U_{\varpi, i}^P)-\chi(t_{i}))^{n(\chi, i)}
\prod_{\chi' \neq \chi} 
(\delta_P^{\frac{1}{2}}(t_{i(\chi')})\pi(U_{\varpi, i(\chi')}^P)-\chi'(t_{i(\chi')}))^{n(\chi')} \phi_K
\]
is a $p$-stabilized vector with respect to $\widehat{P}$,
which satisfies
\[
\pi(U_{\varpi, i}^P)\phi = \delta_P(t_i)^{-\frac{1}{2}}\chi (t_i) \phi,
\]
for  all $i  \in \{1, \ldots, m_P\}$.

\medskip
In the following series of sections, we give examples of $p$-stabilized vectors in various settings. 

\section{$GL_2$-case}
Let $\alpha$ be the simple root of $GL_2$ such that $\alpha:T\lra F^\times,\ {\rm diag}(a,b)\mapsto ab^{-1}$. 
We have a $U_p$-operator 
$U^B_{\w,1}=[It_1I]$, where $t_1=\di(1,\w^{-1})$. 
Let $\pi=\pi(\chi)$ be an unramified principal series representation of ${\rm GL}_2(F)$ where 
$\chi=\chi_1\otimes\chi_2$
and
$\chi_1,\chi_2$ are unramified quasi-characters of $F^\times$.
Then $\pi$ has a non-zero $K$-fixed vector $\phi_0$,
where $K= {\rm GL}_2(\mathfrak{o})$.
We shall give an explicit $p$-stabilized vector for $\pi$.
The semisimplification of $\pi_U$ is 
$$\chi_1\otimes\chi_2+ \chi_2\otimes\chi_1.$$
The element $t_1$ acts on each irreducible component
of $\pi_U$ by $\chi_2(\w^{-1}),\ \chi_1(\w^{-1})$ respectively. 
If $\chi_1(\w^{-1})\neq \chi_2(\w^{-1})$, then 
we have
$\pi_U=(\chi_1\otimes\chi_2)\oplus (\chi_2\otimes\chi_1)$. 
It follows from Proposition~\ref{plus}
and the result in Section 4
that
$$f_1:=(\delta^{\frac{1}{2}}_B(t_1)U^{B}_{\w,1}-\chi_2(\w^{-1}))\phi_0,\ 
f_2:=(\delta^{\frac{1}{2}}_B(t_1)U^{B}_{\w,1}-\chi_1(\w^{-1}))\phi_0$$
are $p$-stabilized vectors with respect to $B$ with the eigenvalues $q^{\frac{1}{2}}\chi_1(\w^{-1}),\ 
q^{\frac{1}{2}}\chi_2(\w^{-1})$
respectively. 

If $\chi_1(\w^{-1})=\chi_2(\w^{-1})$, 
then we have $\chi_1=\chi_2$
since $\chi_1$ and $\chi_2$ are unramified.
In this case, 
any irreducible component of $\pi_U$ is 
isomorphic to
$\chi_1\otimes \chi_1$, but $\pi_U$ is not decomposed as
 $(\chi_1\otimes \chi_1)^{\oplus 2}$. 
This follows from the fact that $\C\simeq {\rm Hom}_{{\rm GL}_2(F)}(\pi,\pi)={\rm Hom}_T(\pi_U,\chi_1\otimes \chi_1)$
(Schur's lemma and Frobenius reciprocity).  
In this case, $f_3:=(\delta^{\frac{1}{2}}_B(t_1)U^{B}_{\w,1}-\chi_1(\w^{-1}))\phi_0$ is a 
$p$-stabilized vector with respect to $B$ with the eigenvalue $q^{\frac{1}{2}}\chi_1(\w^{-1})$. 

We can express $f_i$ in terms of Iwahori fixed vectors as follows. 
Let $\phi$ be a generator of $\pi^K$.  
Choose a basis $\{\phi,\ \phi'=\pi(t^{-1}_1)\phi\}$ of $\pi^I$. 
Then we have 
$$U^B_{\w,1}(\phi,\phi')=(\phi,\phi')
\left(\begin{array}{cc}
a(\phi) & q  \\
-\chi_1(\w^{-1})\chi_2(\w^{-1}) & 0
\end{array}
\right)$$
where $a(\phi)=\chi_1(\w^{-1})+\chi_2(\w^{-1})$. 
Therefore we have $$\chi_1(\w)f_1=\phi-q^{-\frac{1}{2}}\chi_2(\w^{-1})\phi',\ \chi_2(\w)f_2=\phi-q^{-\frac{1}{2}}\chi_1(\w^{-1})\phi'$$
(see Section 9.1 for the relation to the global setting).
It is the same for $f_3$.

\section{$U(2,1)$-case}
Let $U(2,1)$ be the quasi-split unitary group in three variables 
associated to an quadratic extension $E/F$. 
Put $$\Phi=\left(\begin{array}{cccc}
0&0 & 1 \\
0&-1 & 0 \\
1&0 & 0  
\end{array}
\right).$$
We denote by ${}^-$ the conjugate for the non-trivial element in Gal$(E/F)$. 
We realize ${\rm U}(2,1)(F)$ as 
the subgroup of ${\rm GL}_3(E)$ consisting of all $g$ satisfying ${}^t\overline{g}\Phi g=\Phi$. 
Let $B$ be the upper triangular Borel subgroup of $U(2,1)$,
$T$ the diagonal subgroup of $B$
and $K = {\rm U}(2,1)(F)\cap {\rm GL}_3(\mathfrak{o}_E)$,
where $\mathfrak{o}_E$
is the ring of integers in $E$. 
We write $E^1$ for the norm-one subgroup of $E/F$.
Then 
$T$ 
is isomorphic to $E^\times \times E^1$.
Let $\chi$ be an unramified quasi-character of $E^\times$
and let $\one$ denote the trivial character of $E^1$.

Due to \cite{keys},
the corresponding parabolically induced representation
$\mathrm{Ind}_{B(F)}^{{\rm U}(2,1)(F)} (\chi \otimes \one)$
is irreducible except for the following cases:
\begin{enumerate}
\item
$\chi = |\cdot|_E^{\pm}$,
where $|\cdot|_E$ denotes the normalized absolute value
of $E$;

\item
$\chi|_{F^\times} = \omega_{E/F} |\cdot|^{\pm}$,
where $\omega_{E/F}$ is the non-trivial character of $F^\times$
which is trivial on $N_{E/F}(E^\times)$;

\item
$\chi|_{F^\times}$ is trivial 
and $\chi$ is not trivial.
\end{enumerate}
Suppose that
$\pi = \mathrm{Ind}_{B(F)}^{{\rm U}(2,1)(F)} (\chi \otimes \one)$
is irreducible.
Then $\pi$ has a non-zero $K$-fixed vector $\phi_0$.
We shall produce an explicit $p$-stabilized vector for $\pi$.

We fix a uniformizer $\varpi_E$ of $E$
and
set 
\[
t_1 = 
\left(
\begin{array}{ccc}
\varpi_E & 0&0 \\
0& 1 &0 \\
0&0 & \overline{\varpi}_E^{-1}
\end{array}
\right).
\]
Let $U$ be the unipotent radical of $B$.
Then $t_1$ is positive relative to $(B, U)$.
The semisimplification of 
$\pi_U$ is 
$\chi \otimes \one + \overline{\chi}^{-1} \otimes \one$,
where $\overline{\chi}$ denotes the 
quasi-character of $E^\times$ defined by
$\overline{\chi}(x) = \chi(\overline{x})$, for $x \in E^\times$.
The element $t_1$ acts on 
$\chi \otimes \one$ and $\overline{\chi}^{-1} \otimes \one$
by $\chi(\varpi_E)$ 
and $\overline{\chi}^{-1}(\varpi_E)$ respectively.
As in the $GL_2$-case,
by Proposition~\ref{plus} and the result in Section 4,
\[
(\delta_B^{\frac{1}{2}}(t_1)U_{\varpi, 1}^B -\overline{\chi}^{-1}(\varpi_E)) \phi_0,\, 
(\delta_B^{\frac{1}{2}}(t_1)U_{\varpi, 1}^B-\chi(\varpi_E)) \phi_0
\]
are $p$-stabilized vectors
with respect to $B$
with the eigenvalues $q_E\chi(\varpi_E)$ and $q_E\overline{\chi}^{-1}(\varpi_E)$
respectively.

\section{$GSp_4$-case}
Hereafter we follows the notations in \cite{rs}. 
Put $$J=\left(\begin{array}{cccc}
0 &0&0 & 1 \\
0 &0&1 & 0 \\
0 &-1&0 & 0 \\
-1 &0&0 & 0 
\end{array}
\right).$$
We realize ${\rm GSp}_4(F)$ as 
the subgroup of ${\rm GL}_4(F)$ consisting of all $g$ such that 
${}^tgJg=\lambda J$, for some $\lambda\in F^\times$. 
Let $B$ be the Borel subgroup of $GSp_4$ consisting of 
the upper triangular elements, $T$ the diagonal subgroup of $B$
and $U$ the unipotent radical of $B$. 
Let $P$ (respectively $Q$) be the Siegel (respectively Klingen) parabolic subgroup of $GSp_4$
containing $B$. 
Let $P=M_PN_P$ and $Q=M_QN_Q$ be 
Levi decompositions.   
Let $\alpha_1, \alpha_2$ be simple roots  
defined by 
$\alpha_1(t)=ab^{-1}$ and $\alpha_2(t)=b^2c^{-1}$ for $t={\rm diag}(a,b,cb^{-1},ca^{-1})$,
 $a, b, c \in F^\times$. Note that $\Delta_P=\{\alpha_1\}$ and $\Delta_Q=\{\alpha_2\}$. 
Thus there exist two $U_p$-operators $U^B_{\w,i}=[It_iI],\ i=1,2$ where 
$t_1={\rm diag}(1,1,\w^{-1},\w^{-1})$ and $t_2={\rm diag}(1,\w^{-1},\w^{-1},\w^{-2})$. 
Note that $t_1,t_2$ are positive elements in $T_B^-$ 
relative to $(B,I)$, and  
$t_1$ (respectively $t_2$) is a positive element in $T^-_P$ (respectively $T^-_Q$) relative to $(P,K_P)$ (respectively $(Q,K_Q)$).  
In this section,
we give $p$-stabilized vectors 
according to the classification of the
parahori-spherical representations of ${\rm GSp}_4(F)$ by 
Roberts and Schmidt \cite{rs}.

\subsection{Iwahori case}\label{subsec:iwahori}
In this subsection, 
we study $p$-stabilized vectors with respect to 
$\widehat{B} = B$.
The strategy taking here is as follows:

\begin{enumerate}
\item 
Let $(\pi, V)$ be an irreducible admissible representation of ${\rm GSp}_4(F)$ admitting a non-zero Iwahori-fixed vector.
Then $\pi$ is an irreducible constituent of 
some unramified principal series representation 
$\chi_1\times\chi_2\rtimes \sigma$. 
We use the classification of such representations
in 
Table A.15 of \cite{rs}.

\item 
We make up a list of the
simultaneous eigenvalues for $U_p$-operators
$U^B_{\w,1}$, $U^B_{\w,2}$
in terms of the following Satake parameters of $\chi_1\times\chi_2\rtimes \sigma$: 
\begin{eqnarray}\label{eq:satake}
\alpha=\chi_1\chi_2\sigma(\w^{-1}),\ 
\beta=\chi_1\sigma(\w^{-1}),\ 
\gamma=\chi_2\sigma(\w^{-1}),\ 
\delta=\sigma(\w^{-1}).
\end{eqnarray}
By Proposition \ref{plus}, 
the problem is
reduced to the computation of the
simultaneous eigenvalues for $t_1$,
$t_2$ on $(V_U)^{I\cap T}$. 
We note that $V_U = (V_U)^{I\cap T}$.
Table A.3 of \cite{rs} gives the semisimplification
 of the Jacquet module $\pi_{N_P}$ of $\pi$ associated to 
the Siegel parabolic subgroup $P$. 
So we can easily get the semisimplification $\pi_U^{\mathrm{ss}}$
of $\pi_U$
by using the transitivity of Jacquet functors.
The elements $t_1,t_2$ act on 
each irreducible component of $\pi_U^{\mathrm{ss}}$
by its central character.
Thus we can 
obtain the set $S'$ of pairs of  
simultaneous eigenvalues for $(t_1, t_2)$
on $V_U = (V_U)^{I\cap T}$. 
It will turns out that $S'$ is contained in $S$,
where
$S$ is the set of
pairs of  
simultaneous eigenvalues for $(t_1, t_2)$ on 
$(\chi_1\times\chi_2\rtimes \sigma)_U
=
(\chi_1\times\chi_2\rtimes \sigma)_U^{I\cap T}$,
that is,
\[
S= \{(\delta,\delta\gamma),(\delta,\delta\beta),(\alpha,\alpha\beta),
(\alpha,\alpha\gamma),(\gamma,\gamma\delta),(\gamma,\gamma\alpha),(\beta,\beta\delta),(\beta,\beta\alpha)\}.
\]

\item 
Suppose that $S'$ contains just 
$\dim V_U$-elements.
Since $t_1$, $t_2$ generates $T/Z_G$,
this implies that 
$(\pi_U)^{I\cap T}$
is a semisimple and multiplicity-free $T$-module.
We further assume that 
$\pi$ has a non-zero $K$-fixed vector $\phi_K$.
In this case,
given an element $(s, t)$ in $S'$,
the vector 
\[
\phi_{s,t}:=
\prod_{\substack{(s',t')\in S' \\ s' \neq s}}(\delta^{\frac{1}{2}}_B(t_1)U^{B}_{\w,1}-s')
\prod_{\substack{(s',t')\in S' \\ t'\neq t}}
(\delta^{\frac{1}{2}}_B(t_2)U^{B}_{\w,2}-t')\phi_K
\]
is a $p$-stabilized vector with respect to $B$
with the eigenvalues $(\delta^{-\frac{1}{2}}_B(t_1)s, 
\delta^{-\frac{1}{2}}_B(t_2)t)$
because of Proposition~\ref{plus} and the result in Section 4.
We note that Table A. 15 of \cite{rs} gives a list of the
$K$-spherical representations of ${\rm GSp}_4(F)$.
\end{enumerate}

\subsubsection{Case I} 
Let $\chi_1,\chi_2,\sigma$ be unramified quasi-characters of $F^\times$. 
Then the corresponding parabolically induced representation
 $\pi=\chi_1\times\chi_2\rtimes \sigma$ is irreducible 
 if and only if $\chi_1\neq \nu^{\pm 1},\ 
\chi_2\neq \nu^{\pm 1}$ and $\chi_1\neq \nu^{\pm 1}\chi_2^{\pm 1}$. 
Here let us put $\nu=|\cdot|$. 
Due to Table A. 15 of \cite{rs},
$\pi$ has a non-zero $K$-fixed vector.
The semisimplification of $\pi_{N_P}$ is 
given in Table A.3 of \cite{rs}.
We use the notation in section A.3 of \cite{rs}.
Let $r_{T, M_P}$ denote the Jacquet functor 
from the category of the smooth representations of $M_P$
to that of $B$.
Note that the semisimplification of 
$r_{T, M_P}((\chi_1 \times \chi_2)\otimes \sigma)$
is $(\chi_1\otimes\chi_2\otimes\sigma)+
 (\chi_2\otimes\chi_1\otimes\sigma)$.
 Since $r_U = r_{T, M_P} \circ r_{N_P}$,
the semisimplification of $\pi_U$ is
\begin{align*}
(\chi_1\otimes\chi_2\otimes\sigma)+
 (\chi_2\otimes\chi_1\otimes\sigma)+
(\chi^{-1}_1\otimes\chi^{-1}_2\otimes\chi_1\chi_2\sigma)+(\chi^{-1}_2\otimes\chi^{-1}_1\otimes
\chi_1\chi_2\sigma)+\\
(\chi_1\otimes\chi^{-1}_2\otimes\chi_2\sigma)+( \chi^{-1}_2\otimes\chi_1\otimes\chi_2\sigma)+ 
(\chi_2\otimes\chi^{-1}_1\otimes\chi_1\sigma)+(\chi^{-1}_1\otimes\chi_2\otimes
\chi_1\sigma).
\end{align*}
The elements $t_1, t_2$ 
act on each component of the semisimplification of $\pi_U$
by the following pairs of scalars:
$$(\delta,\delta\gamma),(\delta,\delta\beta),(\alpha,\alpha\beta),
(\alpha,\alpha\gamma),(\gamma,\gamma\delta),(\gamma,\gamma\alpha),(\beta,\beta\delta),(\beta,\beta\alpha).$$
The pairs above are pairwise distinct
if and only if $\chi_1 \neq 1$ , $\chi_2 \neq 1$
and $\chi_1 \neq \chi_2^{\pm 1}$.

\subsubsection{Case II}
Let $\chi,\sigma$ be unramified quasi-characters of $F^\times$ such that $\chi^2\neq \nu^{\pm 1}$
and $\chi\neq \nu^{\pm\frac{3}{2}}$. 
Put $\chi_1=\nu^{\frac{1}{2}}\chi$
and $\chi_2=\nu^{-\frac{1}{2}}\chi$. 
In what follows, we consider the irreducible constituents 
$\pi$ of 
$\chi_1\times\chi_2\rtimes \sigma$. 

\medskip
\noindent
{\bf Case IIa.} 
Put $\pi=\chi{\rm St}_{\mathrm{GL}(2)}\rtimes \sigma$.
Noting that $r_{T, M_P}(\chi\mathrm{St}_{\mathrm{GL}(2)}\otimes \sigma) = \nu^{\frac{1}{2}}\chi \otimes
\nu^{-\frac{1}{2}}\chi \otimes \sigma$.
It follows from Table A.3 of \cite{rs} that
the semi-simplification of $\pi_U$ is
$$(\chi_1\otimes\chi_2\otimes\sigma)+(\chi^{-1}_2\otimes\chi^{-1}_1\otimes\chi_1\chi_2\sigma)+
(\chi_1\otimes\chi^{-1}_2\otimes\chi_2\sigma)+(\chi^{-1}_2\otimes\chi_1\otimes\chi_2\sigma).$$
Hence the elements $t_1,t_2$ act on each 
irreducible component of $\pi_U$
by
$$(\delta,\delta\gamma),(\alpha,\alpha\gamma),(\gamma,\gamma\delta),(\gamma,\gamma\alpha)$$  
respectively. 
The pairs above are pairwise distinct
if and only if $\chi^2 \neq 1$.

\medskip
\noindent
{\bf Case IIb. } 
If $\pi=\chi{\mathbf{1}}_{{\rm GL}(2)}\rtimes \sigma$,
then $\pi$ has a non-zero $K$-fixed vector.
Since 
$r_{T, M_P}(\chi\mathbf{1}_{\mathrm{GL}(2)}\otimes \sigma) = \nu^{-\frac{1}{2}}\chi \otimes
\nu^{\frac{1}{2}}\chi \otimes \sigma$,
it follows from Table A.3 of \cite{rs} that
the semisimplification of $\pi_U$ is
$$(\chi_2\otimes\chi_1\otimes\sigma)+
(\chi^{-1}_1\otimes\chi^{-1}_2\otimes\chi_1\chi_2\sigma)
+
(\chi_2\otimes\chi^{-1}_1\otimes\chi_1\sigma)+(\chi^{-1}_1\otimes\chi_2\otimes\chi_1\sigma).$$
The elements $t_1,t_2$ 
act on each component by
$$(\delta,\delta\beta),(\alpha,\alpha\beta),(\beta,\beta\delta),(\beta,\beta\alpha)$$  
respectively. 
The pairs above are pairwise distinct
if and only if $\chi^2 \neq 1$.

\subsubsection{Case III}
Let $\chi$ and $\rho$ be unramified quasi-characters of $F^\times$. 
We assume that $\chi\neq 1$ and $\chi\neq \nu^{\pm 2}$. 
Put $\chi_1=\chi,\ \chi_2=\nu$, and 
$\sigma=\nu^{-\frac{1}{2}}\rho$. 
Next, we consider the irreducible constituents $\pi$ of 
$\chi_1\times\chi_2\rtimes \sigma$. 

\medskip
\noindent
{\bf Case IIIa.} 
If $\pi=\chi\rtimes \rho {\rm St}_{{\rm GSp}(2)}$,
then
the semisimplification of $\pi_U$ is
$$(\chi_1\otimes\chi_2\otimes\sigma)+ (\chi_2\otimes\chi_1\otimes\sigma)+
(\chi_2\otimes\chi^{-1}_1\otimes\chi_1\sigma)+(\chi^{-1}_1\otimes\chi_2\otimes\chi_1\sigma).$$
The actions of $t_i,\ i=1,2 $ on each 
irreducible component of $\pi_U$ are just
$$(\delta,\delta\gamma),(\delta,\delta\beta),(\beta,\beta\delta),(\beta,\beta\alpha)$$  
respectively. 
The pairs above are pairwise distinct
if and only if $\chi \neq \nu^{\pm 1}$.

\medskip
\noindent
{\bf Case IIIb.} 
Suppose that $\pi=\chi\rtimes \rho \mathbf{1}_{{\rm GSp}(2)}$.
Then $\pi$ admits a non-zero
$K$-fixed vector and
the semisimplification of $\pi_U$ is
$$(\chi^{-1}_1\otimes\chi^{-1}_2\otimes\chi_1\chi_2\sigma)
+
(\chi^{-1}_2\otimes\chi^{-1}_1\otimes\chi_1\chi_2\sigma)+
(\chi_1\otimes\chi^{-1}_2\otimes\chi_2\sigma)+(\chi^{-1}_2\otimes\chi_1\otimes\chi_2\sigma).$$
Hence the actions of $t_i,\ i=1,2 $ on each component are 
$$(\alpha,\alpha\beta),(\alpha,\alpha\gamma),(\gamma,\gamma\delta),(\gamma,\gamma\alpha)$$  
respectively. 
The pairs above are pairwise distinct
if and only if $\chi \neq \nu^{\pm 1}$.

\subsubsection{Case IV}
Let $\rho$ be an unramified quasi-character of $F^\times$. 
Put $\chi_1=\nu^2,\ \chi_2=\nu$,
and $\sigma=\nu^{-\frac{3}{2}}\rho$. 
We will consider the irreducible constituents
$\pi$ of 
$\chi_1\times\chi_2\rtimes \sigma$. 
In this case,
$\alpha,\beta,\gamma,\delta$ are different from each other. 
So $\pi_U$ is a semisimple and multiplicity-free
$T$-module.

\medskip
\noindent
{\bf Case IVa.} 
Suppose that $\pi=\rho{\rm St}_{{\rm GSp}(4)}$. 
Then we have 
$\pi_U=\chi_1\otimes\chi_2\otimes\sigma$,
and the elements $t_1,t_2$ act on it
by $(\delta,\delta\gamma)$.
Since $\dim \pi^I = 1$,
any non-zero $I$-fixed vector 
is itself a $p$-stabilized vector with the eigenvalues 
$(\delta_B^{-\frac{1}{2}}(t_1)\delta,
\delta_B^{-\frac{1}{2}}(t_2)\delta\gamma)$.

\medskip
\noindent
{\bf Case IVb.} 
If $\pi=L(\nu^2,\nu^{-1}\rho {\rm St}_{{\rm GSp}(2)})$,
then
$\pi_U$ is isomorphic to
$$(\chi_2\otimes\chi_1\otimes\sigma)\oplus
(\chi_2\otimes\chi^{-1}_1\otimes\chi_1\sigma)\oplus
(\chi^{-1}_1\otimes\chi_2\otimes\chi_1\sigma).$$
Hence $t_i$, $i = 1, 2$ act
on each component by
$$(\delta,\delta\beta), (\beta,\beta\delta), (\beta,\beta\alpha)$$  
respectively.

\medskip
\noindent
{\bf Case IVc.} 
Suppose that $\pi=L(\nu^{\frac{3}{2}}{\rm St}_{{\rm GL}(2)},\nu^{-\frac{3}{2}}\rho)$. 
Then $\pi_U$ is isomorphic to
$$(\chi^{-1}_2\otimes\chi^{-1}_1\otimes\chi_1\chi_2\sigma)\oplus 
(\chi_1\otimes\chi^{-1}_2\otimes\chi_2\sigma)\oplus 
(\chi^{-1}_2\otimes\chi_1\otimes\chi_2\sigma),$$
and the elements $t_1,t_2$
act on each irreducible component by
$$(\alpha,\alpha\gamma), (\gamma,\gamma\delta), (\gamma,\gamma\alpha)$$  
respectively.  

\medskip
\noindent
{\bf Case IVd.} 
If $\pi=\rho \mathbf{1}_{{\rm GSp}(4)}$,
then we have $\dim \pi^I = \dim \pi^K = 1$,
and hence $\pi^I = \pi^K$. We also get $\pi_U=\chi^{-1}_1\otimes\chi^{-1}_2\otimes\chi_1\chi_2\sigma$,
and the elements $t_1,t_2$ act on $\pi_U$ by
$(\alpha,\alpha\beta)$.
Thus any non-zero $K$-fixed vector is a $p$-stabilized vector with the eigenvalues $(\delta_B^{-\frac{1}{2}}(t_1)\alpha,
\delta_B^{-\frac{1}{2}}(t_2)\alpha\beta)$.

\subsubsection{Case V}
Let $\xi$ and $\rho$ be unramified quasi-characters of $F^\times$.
We assume that $\xi^2=1$ and $\xi\neq 1$. 
Put $\chi_1=\nu\xi,\ \chi_2=\xi$,
and $\sigma=\nu^{-\frac{1}{2}}\rho$. 
We consider the irreducible constituents 
$\pi$ of 
$\chi_1\times\chi_2\rtimes \sigma$. 
In this case,
 $\alpha,\beta,\gamma,\delta$ are different from each other. 
This means that $\pi_U$ is a semisimple
and multiplicity-free $T$-module.   

\medskip
\noindent
{\bf Case Va.} 
If $\pi=\delta([\xi,\nu\xi],\nu^{-\frac{1}{2}}\rho)$,
then we have 
$$\pi_U=(\chi_1\otimes\chi_2\otimes\sigma)\oplus  
(\chi_1\otimes\chi^{-1}_2\otimes\chi_2\sigma).$$
The elements $t_1,t_2$ act on each component by
$$(\delta,\delta\gamma), (\gamma,\gamma\delta)$$  
respectively.

\medskip
\noindent
{\bf Case Vb.} 
Suppose that $\pi=L(\nu^{\frac{1}{2}}\xi{\rm St}_{{\rm GL}(2)},\nu^{-\frac{1}{2}}\rho)$. 
Then we have 
$$\pi_U=(\chi^{-1}_2\otimes\chi^{-1}_1\otimes\chi_1\chi_2\sigma)\oplus  
(\chi^{-1}_2\otimes\chi_1\otimes\chi_2\sigma),$$
and hence the elements 
$t_1,t_2$ act on each irreducible component
by
$$(\alpha,\alpha\gamma), (\gamma,\gamma\alpha)$$  
respectively. 

\medskip
\noindent
{\bf Case Vc.} 
In the case when
 $\pi=L(\nu^{\frac{1}{2}}\xi \mathrm{St}_{\mathrm{GL}(2)}, \xi\nu^{-\frac{1}{2}}\rho)$,
we have 
$$\pi_U=(\chi_2\otimes\chi_1\otimes\sigma)\oplus  
(\chi_2\otimes\chi^{-1}_1\otimes\chi_1\sigma).$$
The elements $t_1,t_2$ act on each component
by
$$(\delta,\delta\beta), (\beta,\beta\delta)$$  
respectively. 

\medskip
\noindent
{\bf Case Vd.} 
If
$\pi= L(\nu\xi, \xi\rtimes \nu^{-\frac{1}{2}}\rho)$,
then $\pi$ has a non-zero $K$-fixed vector.
We have 
$$\pi_U=(\chi^{-1}_1\otimes\chi^{-1}_2\otimes\chi_1\chi_2\sigma)\oplus  
(\chi^{-1}_1\otimes\chi_2\otimes\chi_1\sigma).$$
Hence the elements $t_1,t_2$ act on each component by
$$(\alpha,\alpha\beta), (\beta,\beta\alpha)$$  
respectively. 

\subsubsection{Case VI}
Let $\rho$ be an unramified quasi-character of $F^\times$. 
Put $\chi_1 = \nu$, $\chi_2 = \mathbf{1}_{F^\times}$
and $\sigma = \nu^{-1/2}\rho$. 
Finally,
we consider the irreducible constituents $\pi$ of 
$\chi_1\times\chi_2\rtimes \sigma$. 
In this case,
we have $\alpha=\beta,\ \gamma=\delta$ and $\alpha\neq \gamma$. 

\medskip
\noindent
{\bf Case VIa.} 
Suppose that $\pi=\tau(S, \nu^{-\frac{1}{2}}\rho)$. 
Then the semisimplification of $\pi_U$ is
$$(\chi_1\chi_2\otimes\chi_1\otimes\sigma)^{\oplus 2}
+
(\chi_2\otimes\chi_1\otimes\sigma).$$
Hence the actions of $t_i$,
$i= 1,2$ on each component are 
$$(\gamma,\gamma^2), (\gamma,\gamma^2), (\gamma,\gamma\alpha)$$  
respectively.  

\medskip
\noindent
{\bf Case VIb.} 
If $\pi= \tau(T, \nu^{-\frac{1}{2}}\rho)$,
then we have $\dim \pi^{K_P} = \dim \pi^I = 1$
by Table A. 15 in \cite{rs},
and hence $\pi^I = \pi^{K_P}$.
We also get 
$\pi_U=\chi_2\otimes\chi_1\otimes\sigma$,
and the elements $t_1, t_2$ act on $\pi_U$ by
$(\gamma,\gamma\alpha)$.
Therefore any non-zero $K_P$-spherical vector
is a  $p$-stabilized vector with the eigenvalues 
$(\delta_B^{-\frac{1}{2}}(t_1)\gamma,
\delta_B^{-\frac{1}{2}}(t_2)\gamma\alpha)$.

\medskip
\noindent
{\bf Case VIc.} 
Suppose that
 $\pi= L(\nu^{\frac{1}{2}}\mathrm{St}_{\mathrm{GL}(2)}, \nu^{-\frac{1}{2}} \rho)$.  
Then we get
$\pi_U=\chi_2\otimes\chi^{-1}_1\otimes\chi_1\sigma$,
and $t_i$, $i = 1,2$ acts on it by
$(\alpha,\gamma\alpha)$.
By Table A. 15 in \cite{rs},
the space $\pi^{K_Q}$ is one-dimensional. 
Therefore
a non-zero $K_Q$-fixed vector 
is  a $p$-stabilized vector with the eigenvalues 
$(\delta_B^{-\frac{1}{2}}(t_1)\alpha,
\delta_B^{-\frac{1}{2}}(t_2)\gamma\alpha)$. 

\medskip
\noindent
{\bf Case VId.}
If
$\pi=L(\nu, \mathbf{1}_{F^\times}\rtimes \nu^{-1/2}\rho)$,
then $\pi$ has a non-zero $K$-fixed vector.
The semisimplification of $\pi_U$ is
$$(\chi^{-1}_1\otimes\chi_2\otimes\chi_1\sigma)^{\oplus 2}
+
(\chi_2\otimes\chi^{-1}_1\otimes\chi_1\sigma).$$
Hence the actions of $t_i,\ i=1,2 $ on each component are 
just
$$(\alpha,\alpha^2), (\alpha,\alpha^2), (\alpha,\gamma\alpha)$$  
respectively.  

\medskip
In Table 1,
we list the simultaneous eigenvalues 
for $t_1$, $t_2$ on $\pi_U$,
where $\pi$ is a representation in groups I-V.
For a representation in group VI,
we always have
$\alpha=\beta,\ \gamma=\delta$ and $\alpha\neq \gamma$. 
So for such representations,
we list the multiplicity of 
the simultaneous eigenvalues 
for $t_1$, $t_2$ on $\pi_U$
in Table 2.

\begin{table}[htbp]
\label{tab1}
\begin{center}
\begin{tabular}{|c|cccccccc|c|c|}
\hline
rep. & $(\alpha, \alpha \beta)$ &
$(\alpha, \alpha\gamma)$ & 
$(\beta, \beta \alpha)$ & $(\beta, \beta \delta)$
& $(\gamma, \gamma\alpha)$ & $(\gamma, \gamma\delta)$
& $(\delta, \delta \beta)$ & $(\delta, \delta\gamma)$
& $\dim \pi^I$ & $\dim \pi^K$\\
\hline
I & $\bigcirc$ & $\bigcirc$ &  $\bigcirc$ & $\bigcirc$ &
$\bigcirc$ & $\bigcirc$ &  $\bigcirc$ & $\bigcirc$ & 8 & 1\\
\hline
IIa & - & $\bigcirc$   & - & - & $\bigcirc$   & $\bigcirc$   & - & $\bigcirc$  & 4& 0\\
IIb  & $\bigcirc$ & - &$\bigcirc$ & $\bigcirc$ & - & - & $\bigcirc$ & - & 4& 1\\
\hline
IIIa   & - & - & $\bigcirc$ & $\bigcirc$ & - & - & $\bigcirc$ & $\bigcirc$ & 4& 0\\
IIIb    & $\bigcirc$ & $\bigcirc$ & - & - & $\bigcirc$ & $\bigcirc$ & - & - & 4& 1\\
\hline
IVa & - & - & - & - & - & - & - & $\bigcirc$ & 1 & 0\\
IVb & - & - & $\bigcirc$ & $\bigcirc$ & - & - & $\bigcirc$ & - & 3
& 0\\
IVc  & - & $\bigcirc$ & - & - & $\bigcirc$ & $\bigcirc$ & - & - & 3
& 0\\
IVd  & $\bigcirc$ & - & - & - & - & - & - & - & 1 & 1\\
\hline
Va & - & - & - & - & - & $\bigcirc$ & - & $\bigcirc$ & 2& 0\\
Vb  & - & $\bigcirc$ & - & - &$\bigcirc$ & - & - & - & 2& 0\\
Vc   & - & - & - & $\bigcirc$ & - & - & $\bigcirc$ & - & 2& 0\\
Vd    & $\bigcirc$ & - & $\bigcirc$ & - & - & - & - & - &2& 1\\
\hline
\end{tabular}
\end{center}
\caption{}
\end{table}

\begin{table}[htbp]
\label{tab2}
\begin{center}
\begin{tabular}{|c|cccc|c|c|}
\hline
representation & $(\alpha, \alpha^2)$ & $(\alpha, \alpha\gamma)$ & $(\gamma, \gamma \alpha)$ & $(\gamma, \gamma^2)$ & $\dim \pi^I$ & $\dim \pi^K$\\
\hline
VIa & 0 &0 &  1& 2 & 3 & 0\\
VIb & 0 &  0&  1 & 0 & 1& 0\\
VIc & 0 & 1 & 0 & 0 & 1& 0\\
VId & 2 & 1 & 0 & 0 & 3& 1\\
\hline
\end{tabular}
\end{center}
\caption{}
\end{table}

\subsection{Siegel parahoric case}
In this subsection, we 
compute the eigenvalues of 
the $U_p$-operator $U_{\w, 1}^P$,
where $P$ denotes the Siegel parabolic subgroup
of ${\rm GSp}_4(F)$.
The strategy is as follows:

\begin{enumerate}
\item  
We use the classification of 
the irreducible smooth representations $(\pi, V)$ of ${\rm GSp}_4(F)$
admitting $K_P$-fixed vectors in 
Table A.15 of \cite{rs}.
As in the previous subsection,
we realize $\pi$ as an irreducible constituent of  
some unramified principal series representation  $\chi_1\times \chi_2 \rtimes \sigma$. 
We use the same notation as in subsection~\ref{subsec:iwahori}.

\item
We compute the set $S'$ of
eigenvalues of $t_1$ on 
$(V_{N_P})^{M_P\cap K_P}$
as follows:
The semisimplification $\pi_{N_P}^{\mathrm{ss}}$ of 
the Jacquet module $\pi_{N_P}$ of $\pi$
associated to $P$ is given in Table A.3 of \cite{rs}. 
We denote by $\pi_{N_P,M_P\cap K_P}$
the $M_P$-submodule of $\pi_{N_P}^{\mathrm{ss}}$
spanned by the $M_P\cap K_P$-fixed vectors
in $\pi_{N_P}^{\mathrm{ss}}$.
Note that for any irreducible admissible representation 
$\tau$ of $M_P \simeq {\rm GL}_2(F) \times F^\times$,
we have $\dim \tau^{M_P\cap K_P} \leq 1$.
Thus the length of $\pi_{N_P,M_P\cap K_P}$
is equal to $\dim (V_{N_P})^{M_P\cap K_P}$.
Since the element $t_1$ lies in the center of $M_P$,
the eigenvalues of $t_1$ on $(V_{N_P})^{M_P\cap K_P}$
is just those of $t_1$ on  $\pi_{N_P,M_P\cap K_P}$.
So 
we can easily compute the eigenvalues of $t_1$ on 
$(V_{N_P})^{M_P\cap K_P}$
because $t_1$ 
acts on each irreducible component of 
$\pi_{N_P,M_P\cap K_P}$
by the central character.
It will turns out that $S'$ is contained in 
$S=\{\alpha, \beta, \gamma, \delta\}$,
where $\alpha, \beta, \gamma, \delta$
are the Satake parameters of $\chi_1\times \chi_2 \rtimes \sigma$ defined in (\ref{eq:satake}).

\item 
Suppose that $S'$ contains just 
$\dim (\pi_{N_P})^{K_P\cap M_P}$-elements.
Then
$(\pi_{N_P})^{K_P\cap M_P}$
is a semisimple and multiplicity-free $Z_{M_P}$-module
because $t_1$ generates $Z_{M_P}$ modulo $Z_{M_P}(\mathfrak{o})$.
We further assume that 
$\pi$ has a non-zero $K$-fixed vector $\phi_K$.
Then,
given an element $s$ in $S'$,
the vector 
\[
\phi_{s}:=
\prod_{\substack{s'\in S' \\ s' \neq s}}(\delta^{\frac{1}{2}}_P(t_1)U^{P}_{\w,1}-s')
\phi_K
\]
is a $p$-stabilized vector with respect to $\widehat{P}=Q$
with the eigenvalue $\delta^{-\frac{1}{2}}_P(t_1)s$
because of Proposition~\ref{plus} and the result in Section 4.
\end{enumerate}

\subsubsection{Case I}
In this case, $\pi_{N_P,M_P\cap K_P}$ is 
\[
(\chi_1 \times \chi_2) \otimes \sigma
+
(\chi_1^{-1} \times \chi_2^{-1}) \otimes \chi_1\chi_2\sigma
+(\chi_1 \times \chi_2^{-1}) \otimes \chi_2\sigma
+(\chi_2 \times \chi_1^{-1}) \otimes \chi_1\sigma.
\]
The element  $t_1$ acts on each component 
of $\pi_{N_P,M_P\cap K_P}$ by
\[
\delta, \alpha, \gamma, \beta
\]
respectively. 

\subsubsection{Case II}

\medskip
\noindent
{\bf Case IIa}.
In this case, we have 
$\pi_{N_P,M_P\cap K_P}=
(\chi_1 \times \chi_2^{-1}) \otimes \chi_2\sigma
$. 
The element $t_1$ acts on it by
$\gamma$.
Since $\dim V^{K_P} = 1$,
any non-zero $K_P$-fixed vector is  a $p$-stabilized vector 
with respect to $\widehat{P}$.

\medskip
\noindent
{\bf Case IIb}.
In this case, we have 
\[\pi_{N_P,M_P\cap K_P}=
\chi \mathbf{1}_{\mathrm{GL}(2)} \otimes \sigma
+\chi^{-1} \mathbf{1}_{\mathrm{GL}(2)} \otimes \chi^2 \sigma
+(\chi_2 \times \chi_1^{-1}) \otimes \chi_1\sigma,
\]
and $t_1$ acts on each component by
$
\delta, \alpha, \beta$
respectively. 
\subsubsection{Case III}

\medskip
\noindent
{\bf Case IIIa}.
In this case, we obtain
\[\pi_{N_P,M_P\cap K_P}=
(\chi_1 \times \chi_2) \otimes \sigma
+(\chi_2 \times \chi_1^{-1}) \otimes \chi_1\sigma,
\]
and  $t_1$ acts on each component by
$\delta$ and $\beta$ respectively. 

\medskip
\noindent
{\bf Case IIIb}.
We get
\[\pi_{N_P,M_P\cap K_P}=
(\chi_2^{-1} \times \chi_1^{-1}) \otimes \chi_1\chi_2\sigma
+(\chi_1 \times \chi_2^{-1}) \otimes \chi_2\sigma.
\]
The element  $t_1$ acts on each component by
$\alpha$ and $\gamma$ 
respectively. 
\subsubsection{Case IV}

\medskip
\noindent
{\bf Case IVa}.
The representations in case IVa
have no $K_P$-fixed vectors.

\medskip
\noindent
{\bf Case IVb}.
We obtain
\[\pi_{N_P,M_P\cap K_P}=
\nu^{3/2}\mathbf{1}_{\mathrm{GL}(2)} \otimes \sigma
+(\chi_2 \times \chi_1^{-1}) \otimes \chi_1\sigma. 
\]
The element  $t_1$ acts on each component 
by $\delta$ and $\beta$ respectively. 

\medskip
\noindent
{\bf Case IVc.}
In this case, 
$t_1$ acts on 
$\pi_{N_P,M_P\cap K_P}=
(\chi_1 \times \chi_2^{-1}) \otimes \chi_2\sigma$
by $\gamma$. 

\medskip
\noindent
{\bf Case IVd}.
The element $t_1$ acts on 
$\pi_{N_P,M_P\cap K_P}=
\nu^{-3/2}\mathbf{1}_{\mathrm{GL}(2)} \otimes \chi_1\chi_2 \sigma$  by $\alpha$.

\subsubsection{Case V}

\medskip
\noindent
{\bf Case Va}. 
In this case,
$\pi$ admits no $K_P$-fixed vectors.

\medskip
\noindent
{\bf Case Vb}.
We have 
$\pi_{N_P,M_P\cap K_P}=
\nu^{1/2}\xi \mathbf{1}_{\mathrm{GL}(2)} \otimes \chi_2 \sigma$,
and $t_1$ acts on it by
$\gamma$.

\medskip
\noindent
{\bf Case Vc}.
The element $t_1$ acts on 
$\pi_{N_P,M_P\cap K_P}=
\nu^{1/2}\xi \mathbf{1}_{\mathrm{GL}(2)} \otimes \sigma$
by $\delta$.

\medskip
\noindent
{\bf Case Vd}.
In this case, we obtain
\[
\pi_{N_P,M_P\cap K_P}
=
\nu^{-1/2}\xi \mathbf{1}_{\mathrm{GL}(2)} \otimes \chi_1 \sigma
\oplus
\nu^{-1/2}\xi \mathbf{1}_{\mathrm{GL}(2)} \otimes \chi_1 \chi_2 \sigma.
\]
So  $t_1$ acts on each component by
$\beta$ and $\alpha$ 
respectively.

\subsubsection{Case VI}

\medskip
\noindent
{\bf Case VIa}.
The element $t_1$ acts on 
$\pi_{N_P,M_P\cap K_P}=
\nu^{1/2} \mathbf{1}_{\mathrm{GL}(2)} \otimes \sigma$ 
by $\gamma$.

\medskip
\noindent
{\bf Case VIb}.
We have 
$\pi_{N_P,M_P\cap K_P}=
\nu^{1/2} \mathbf{1}_{\mathrm{GL}(2)} \otimes \sigma$.
The element $t_1$ acts on it by
$\gamma$.

\medskip
\noindent
{\bf Case VIc}.
In this case,
$\pi$ has no $K_P$-fixed vectors.

\medskip
\noindent
{\bf Case VId}.
We get $\pi_{N_P,M_P\cap K_P}=
(\nu^{-1/2} \mathbf{1}_{\mathrm{GL}(2)} \otimes \chi_1 \sigma)^{\oplus 2}$
and  $t_1$ acts on each component by
$\alpha$.

\medskip
In Table 3,
we list the eigenvalues 
for $t_1$ on $(\pi_{N_P})^{M_P\cap K_P}$,
where $\pi$ is a representation in groups I-V.
For representations in group VI,
we list the multiplicity of 
the eigenvalues 
for $t_1$ on $(\pi_{N_P})^{M_P\cap K_P}$
in Table 4.

\begin{table}[htbp]
\label{tab3}
\begin{center}
\begin{tabular}{|c|cccc|c|c|}
\hline
representation & $\alpha$ &  $\beta$ & $\gamma$ & $\delta$&
$\dim \pi^{K_P}$ & $\dim \pi^K$\\
\hline
I & $\bigcirc$ & $\bigcirc$ &  $\bigcirc$ & $\bigcirc$ & 4 & 1\\
\hline
IIa & - & - &  $\bigcirc$ & - & 1 & 0\\
IIb & $\bigcirc$ & $\bigcirc$ &  - & $\bigcirc$ & 3 & 1\\
\hline
IIIa & - & $\bigcirc$ &  - & $\bigcirc$ & 2 & 0\\
IIIb & $\bigcirc$ &- &  $\bigcirc$ & - & 2 & 1\\
\hline
IVa & - & - &  - & - & 0 & 0\\
IVb & - & $\bigcirc$ &  - & $\bigcirc$ & 2 & 0\\
IVc & - & - &  $\bigcirc$ & - & 1 & 0\\
IVd & $\bigcirc$ & - &  - & - & 1 & 1\\
\hline
Va & - & - &  - & - & 0 & 0\\
Vb & - & - &  $\bigcirc$ & - & 1 & 0\\
Vc & - & - &  - & $\bigcirc$ & 1 & 0\\
Vd & $\bigcirc$ & $\bigcirc$ &  - & - & 2 & 1\\
\hline
\end{tabular}
\end{center}
\caption{}
\end{table}

\begin{table}[htbp]
\label{tab4}
\begin{center}
\begin{tabular}{|c|cc|c|c|}
\hline
representation & $\alpha$ & $\gamma$ & $\dim \pi^{K_P}$ 
& $\dim \pi^K$\\
\hline
VIa & 0 & 1 & 1 & 0\\
VIb & 0 & 1 & 1 & 0\\
VIc & 0 & 0 & 0 & 0\\
VId & 2 & 0 & 2 & 1\\
\hline
\end{tabular}
\end{center}
\caption{}
\end{table}

\subsection{Klingen parahoric case}
In this subsection, we 
compute the eigenvalues of 
the $U_p$-operator $U^Q_{\w, 1}
=[K_{Q}t_2K_Q]$.
The strategy is exactly same 
to that for the Siegel parahoric case.
So we shall be brief here.

\begin{enumerate}
\item  
According to the classification
in
Table A.15 of \cite{rs},
we 
realize the representations $(\pi, V)$ of ${\rm GSp}_4(F)$ 
which admit $K_Q$-fixed vectors
as an irreducible constituent of  
some unramified principal series representation  $\chi_1\times \chi_2 \rtimes \sigma$.

\item
We compute the set $S'$ of
eigenvalues of $t_2$ on 
$(V_{N_Q})^{M_Q\cap K_Q}$
as follows:
The semisimplification $\pi_{N_Q}^{\mathrm{ss}}$ of 
$\pi_{N_Q}$  is given in Table A.3 of \cite{rs}. 
We denote by $\pi_{N_Q,M_Q\cap K_Q}$
the $M_Q$-submodule of $\pi_{N_Q}^{\mathrm{ss}}$
generated by the $M_Q\cap K_Q$-fixed vectors.
Note that for any irreducible admissible representation 
$\tau$ of $M_Q \simeq F^\times \times {\rm GSp}_2(F)$,
we have $\dim \tau^{M_Q\cap K_Q} \leq 1$.
So the length of $\pi_{N_Q,M_Q\cap K_Q}$
is equal to $\dim (V_{N_Q})^{M_Q\cap K_Q}$.
Since $t_2 \in Z_{M_Q}$,
the eigenvalues of $t_2$ on $(V_{N_Q})^{M_Q\cap K_Q}$
is just those of $t_2$ on  $\pi_{N_Q,M_Q\cap K_Q}$.
So 
we can compute the eigenvalues of $t_2$ on 
$(V_{N_Q})^{M_Q\cap K_Q}$
because $t_2$ 
acts on each irreducible component of 
$\pi_{N_Q,M_Q\cap K_Q}$
by the central character.
It will turns out that $S'$ is contained in 
$S=\{\alpha \beta, \alpha\gamma, \delta \beta, \delta\gamma\}$,
where $\alpha, \beta, \gamma, \delta$
are the Satake parameters of $\chi_1\times \chi_2 \rtimes \sigma$ defined in (\ref{eq:satake}).

\item 
If $S'$ contains just 
$\dim (\pi_{N_Q})^{K_Q\cap M_Q}$-elements,
then
$(\pi_{N_Q})^{K_Q\cap M_Q}$
is a semisimple and multiplicity-free $Z_{M_Q}$-module.
We further assume that 
$\pi$ has a non-zero $K$-fixed vector $\phi_K$.
Then,
given an element $s$ in $S'$,
the vector 
\[
\phi_{s}:=
\prod_{\substack{s'\in S' \\ s' \neq s}}(\delta^{\frac{1}{2}}_Q(t_2)U^{Q}_{\w,1}-s')
\phi_K
\]
is a $p$-stabilized vector with respect to $\widehat{Q}=P$
with the eigenvalue $\delta^{-\frac{1}{2}}_Q(t_2)s$
because of Proposition~\ref{plus} and the result in Section 4.
\end{enumerate}

\subsubsection{Case I}
In this case, we have 
\[\pi_{N_Q,M_Q\cap K_Q}=
\chi_1 \otimes (\chi_2 \rtimes \sigma)
+
\chi_2 \otimes (\chi_1 \rtimes \sigma)
+
\chi_2^{-1} \otimes (\chi_1 \rtimes \chi_2 \sigma)
+
\chi_1^{-1} \otimes (\chi_2 \rtimes \chi_1 \sigma).
\]
The element   $t_2$ 
acts on each irreducible component of $\pi_{N_Q,M_Q\cap K_Q}$ by
\[
\delta \gamma,  \beta \delta, \alpha \gamma,
\alpha \beta
\]
respectively. 

\subsubsection{Case II}
\medskip
\noindent
{\bf Case IIa}.
In this case, we get
$\pi_{N_Q,M_Q\cap K_Q}=
\chi_1 \otimes (\chi_2 \rtimes \sigma)
+
\chi_2^{-1} \otimes (\chi_1 \rtimes \chi_2 \sigma)$.
The element  $t_2$ acts on each component  by
$\delta \gamma$ and $\alpha \gamma$
respectively.

\medskip
\noindent
{\bf Case IIb}.
The element $t_2$ acts on 
each component of 
$\pi_{N_Q,M_Q\cap K_Q}=\chi_2 \otimes (\chi_1 \rtimes \sigma)
+
\chi_1^{-1} \otimes (\chi_2 \rtimes \chi_1 \sigma)
$
by 
$\beta \delta$ and $\alpha \beta$ respectively.

\subsubsection{Case III}

\medskip
\noindent
{\bf Case IIIa}.
In this case, $t_2$ acts on
$\pi_{N_Q,M_Q\cap K_Q}=
\chi_2 \otimes (\chi_1 \rtimes \sigma)$
by 
$\beta \delta $.

\medskip
\noindent
{\bf Case IIIb}.
We have 
$\pi_{N_Q,M_Q\cap K_Q}=
\chi \otimes \rho \mathbf{1}_{\mathrm{GSp}(2)}
+
\chi^{-1} \otimes \chi \rho \mathbf{1}_{\mathrm{GSp}(2)}
+
\chi_2^{-1} \otimes (\chi_1 \rtimes \chi_2 \sigma)
$,
and  $t_2$ acts on each component of it
by 
$\delta \gamma,  \alpha\beta$ and $\alpha \gamma$ 
respectively.

\subsubsection{Case IV}

\medskip
\noindent
{\bf Case IVa}.
In this case,
$\pi$ has no $K_Q$-fixed vectors.

\medskip
\noindent
{\bf Case IVb}.
The element $t_2$ acts on 
$\pi_{N_Q,M_Q\cap K_Q}=
\chi_2 \otimes (\chi_1 \rtimes \sigma)$
by
$\beta \delta$.

\medskip
\noindent
{\bf Case IVc}.
We have 
$\pi_{N_Q,M_Q\cap K_Q}=
\nu^2\otimes \nu^{-1}\rho
\mathbf{1}_{\mathrm{GSp}(2)}
+
\chi_2^{-1} \otimes (\chi_1 \rtimes \chi_2 \sigma)
$,
and $t_2$ acts on each component by
$\delta \gamma$ and $\alpha \gamma$ 
respectively.

\medskip
\noindent
{\bf Case IVd}.
The element $t_2$ acts on 
$\pi_{N_Q,M_Q\cap K_Q}=\nu^{-2}\otimes \nu\rho
\mathbf{1}_{\mathrm{GSp}(2)}$
by
$\alpha\beta$.

\subsubsection{Case V}

\medskip
\noindent
{\bf Case Va}.
In this case, $t_2$ acts on 
$\pi_{N_Q,M_Q\cap K_Q}=
\chi_1 \otimes (\chi_2 \rtimes \sigma)$
by
$\delta \gamma$.

\medskip
\noindent
{\bf Case Vb}.
In this case, $t_2$ acts on 
$\pi_{N_Q,M_Q\cap K_Q}=
\chi_2^{-1} \otimes (\chi_1 \rtimes \chi_2 \sigma)$
by
$\alpha \gamma$.

\medskip
\noindent
{\bf Case Vc}.
The element $t_2$ acts on 
$\pi_{N_Q,M_Q\cap K_Q}=
\chi_2 \otimes (\chi_1 \rtimes \sigma)$
by
$\beta \delta$. 

\medskip
\noindent
{\bf Case Vd}.
In this case, 
$t_2$ acts on 
$\pi_{N_Q,M_Q\cap K_Q}=
\nu^{-1/2}\xi \otimes (\xi \rtimes \nu^{1/2} \rho)$
by
$\alpha \beta$.

\subsubsection{Case VI}

\medskip
\noindent
{\bf Case VIa}.
We have
$\pi_{N_Q,M_Q\cap K_Q}=
\chi_1 \otimes (\chi_2 \rtimes \sigma)$,
and hence $t_2$ acts on it by
$\gamma^2$.

\medskip
\noindent
{\bf Case VIb}.
In this case,
$\pi$ has no $K_Q$-fixed vectors.

\medskip
\noindent
{\bf Case VIc}.
The element $t_2$ acts on 
$\pi_{N_Q,M_Q\cap K_Q}=
\mathbf{1}_{F^\times} \otimes \rho \mathbf{1}_{{\rm GSp}(2)}$ 
by $\alpha \gamma$.

\medskip
\noindent
{\bf Case VId}.
In this case, $t_2$ acts on 
each component of
$\pi_{N_Q,M_Q\cap K_Q}=
\mathbf{1}_{F^\times} \otimes \rho \mathbf{1}_{{\rm GSp}(2)}
+
\chi_1^{-1} \otimes (\chi_2 \rtimes \chi_1 \sigma)
$
by
$\alpha \gamma$ and $\alpha^2$ respectively. 

\medskip
In Table 5,
we list the eigenvalues 
for $t_2$ on $(\pi_{N_Q})^{M_Q\cap K_Q}$,
where $\pi$ is a representation in groups I-V.
For representations in group VI,
we list the multiplicity of 
the eigenvalues 
for $t_2$ on $(\pi_{N_Q})^{M_Q\cap K_Q}$
in Table 6.

\begin{table}[htbp]
\label{tab5}
\begin{center}
\begin{tabular}{|c|cccc|c|c|}
\hline
representation & $\alpha\beta$ &  $\alpha\gamma$ & $\delta \beta$ & $\delta\gamma$&
$\dim \pi^{K_Q}$ & $\dim \pi^K$\\
\hline
I & $\bigcirc$ & $\bigcirc$ &  $\bigcirc$ & $\bigcirc$ & 4 & 1\\
\hline
IIa & - & $\bigcirc$ &  - & $\bigcirc$ & 2 & 0\\
IIb & $\bigcirc$ & - &  $\bigcirc$ & - & 2 & 1\\
\hline
IIIa & - & - &  $\bigcirc$  & -& 1 & 0\\
IIIb & $\bigcirc$ &$\bigcirc$  &  - & $\bigcirc$  & 3 & 1\\
\hline
IVa & - & - &  - & - & 0 & 0\\
IVb & - & - & $\bigcirc$ &  - & 1 & 0\\
IVc & - & $\bigcirc$  &  -& $\bigcirc$  & 2 & 0\\
IVd & $\bigcirc$ & - &  - & - & 1 & 1\\
\hline
Va & - & - &  - & $\bigcirc$ & 1 & 0\\
Vb & - &  $\bigcirc$ & - & - & 1 & 0\\
Vc & - & - &   $\bigcirc$ & -& 1 & 0\\
Vd & $\bigcirc$ & - &  - & - & 1 & 1\\
\hline
\end{tabular}
\end{center}
\caption{}
\end{table}

\begin{table}[htbp]
\label{tab6}
\begin{center}
\begin{tabular}{|c|ccc|c|c|}
\hline
representation & $\alpha^2$ & $\alpha \gamma $&$\gamma^2$ & $\dim \pi^{K_Q}$ & $\dim \pi^K$\\
\hline
VIa & 0 & 0 & 1& 1 & 0\\
VIb & 0 & 0 & 0& 0 & 0\\
VIc & 0 & 1 & 0& 1 & 0\\
VId & 1 & 1 & 0& 2 & 1\\
\hline
\end{tabular}
\end{center}
\caption{}
\end{table}

\section{A relation to global objects}
In this section, we will be concerned with global objects and give an 
answer why we consider the actions of the positive elements in local settings. 
Basic references of this section are  \cite{Milne} and \cite{Goresky} (see also Section 5 of II in \cite{Satake}). 

Let $G$ be a non-compact connected semisimple algebraic 
group over $\Q$ and $K$ the maximal compact subgroup of $G(\R)$. 
Let $K_\C$ be the complexification of $K$. 
Put $D=G(\R)/K$.   
Let $\lambda:K_\C\lra {\rm Aut}(E_\lambda)$ be an algebraic representation 
on a complex vector space $E_\lambda$. 
Then we obtain a homogeneous vector bundle $\mathbb{E}_\lambda=G(\R)\times_{K,\lambda|_{K}}E_\lambda ={\rm Ind}^{G(\R)}_K \lambda|_{K}$ on $D$. 
Since $D$ is simply connected, there is a (smooth) trivialization 
$$\mathbb{E}_\lambda\stackrel{\sim}{\lra} D\times E_\lambda.$$
With respect to this trivialization,
the action of $\gamma\in G(\R)$ on $D\times E_\lambda$ is given by 
$$\gamma(z,v)=(\gamma z, J_\lambda(\gamma,z)v)$$
for some smooth mapping $J_\lambda:G\times D\lra {\rm Aut}(E_\lambda)$ which is so called an automorphic factor.  
The action of $G(\R)$ induces the cocycle condition 
$J_\lambda(\gamma_1\gamma_2,z)=J_\lambda(\gamma_1,\gamma_2z)J_\lambda(\gamma_2,z)$ for all $\gamma_1,\gamma_2\in G(\R)$.  

For a holomorphic $E_\lambda$-valued function $f$ on $D$, the action of $g \in G(\R)$ is defined by 
$$(g^{-1}f)(z)=J^{-1}_\lambda(g,z)f(gz).$$ 
For a discrete subgroup $\Gamma$ of $G(\Q)$ and a finite character 
$\chi:\Gamma\lra \C^\times$,  we say that a holomorphic function $f:D\lra E_\lambda$ is 
an automorphic form of weight $\lambda$ and level $\Gamma$ 
with the character $\chi$  if it satisfies $\gamma f=\chi(\gamma)f$ for any $\gamma\in\Gamma$. 
We denote by $M_\lambda(\Gamma,\chi)$ the space of automorphic forms of weight $\lambda$ and level $\Gamma$ 
with the character $\chi$. 

Next we construct adelic forms from classical automorphic forms. 
Let $G$ be a connected reductive group over $\Q$ and $Z_G$ the center of $G$. Then its derived group 
$G^{{\rm der}}$ is a semisimple connected algebraic group. 
Consider the following exact sequence 
$$1\lra G^{{\rm der}}\lra G\stackrel{\nu}{\lra} T:=G/G^{{\rm der}}\lra 1$$
where $T$ is a torus (cf. p.303 of \cite{Milne}) . 
If $G=GL_2$ (respectively $GSp_{2n}$), then $T=\mathbb{G}_m$ and $\nu$ is the determinant map 
(respectively the similitude character). 
Let $G(\R)^+$ be the connected component of the identity element in $G(\R)$ with respect to the real topology. 
For simplicity, in what follows we assume that $G(\R)^+=Z_G(\R)^+G^{{\rm der}}(\R)^+$ 
(for example, $G=GSp_{2n}, GL_2, U(p,q), p,q>0$ satisfy this condition). 

Let $\A$ be the adele ring of $\Q$ and $\A_f$ the finite adele of $\Q$.  
Let $K$ be a compact open subgroup of $G(\A_f)$. 
Assume that $\nu(K)\supset T(\hat{\Z})$. 
Then it follows from the strong approximation theorem for $G^{{\rm der}}$ that 
\begin{equation}
\label{sat} G(\A)=G(\Q)G(\R)^+K=G(\Q)Z_G(\R)^+G^{{\rm der}}(\R)^+K.
\end{equation}
Let  $K^{(1)}_\infty$ be the maximal compact subgroup of $G^{{\rm der}}$ and 
set $K_\infty=Z(\R)^+K^{(1)}_\infty$. 
Then $D:=G(\R)^+/K_\infty=G^{{\rm der}}(\R)/K^{(1)}_\infty$ is the bounded symmetric domain 
endowed with an involution $\iota$.   
Let $I\in D$ be the fixed point of $\iota$. From the description of $D$ as above, 
$Z(\R)^+$ acts on $I$ trivially. 
Put $\Gamma=G^{{\rm der}}(\Q)\cap K$. 
For an automorphic form $f\in M_\lambda(\Gamma,\chi)$, 
we define the function $F_f:G(\A)\lra E_\lambda$ as follows. 
By (\ref{sat}), it is possible  to write a given $g\in G(\A)$ as 
$g=az_\infty g_\infty k$ with $a\in G(\Q)$, $z_\infty\in Z_G(\R)^+$  , $g_\infty\in G^{{\rm der}}(\R)^+$, and $k\in K$. 
Then we put 
$$F_f(g)=J^{-1}_\lambda(g_\infty,I)f(g_\infty I).$$
We can check that $F_f$ is an automorphic form on $G(\A)$ in 
the sense of \cite{BJ}. 
Let $\pi_f$ be the maximal irreducible subquotient of the representation of $G(\A)$ generated by $F_f$. 

We now consider Hecke operators and then compare 
them with those in local settings. 
For $\alpha \in G(\Q)^+=G(\Q)\cap G(\R)^+$, consider the double coset 
$$T(\alpha):=\Gamma \alpha \Gamma=\bigcup_i\Gamma \alpha_i.$$
We write $\alpha_i=\alpha_{i,0} \alpha_{i,1}$, 
where $\alpha_{i,0}\in Z(\R)^+$ and $\alpha_{i,1}\in 
G^{{\rm der}}(\R)^+$. 
Then we define the action of $T(\alpha)$ on $f$ by 
$$T(\alpha)f(z)=n(\alpha)_\lambda \sum_i (\alpha_{i,1} f)(z)$$
where 
$n(\alpha)_\lambda\in \Q^\times $ is a normalized factor depending on $\lambda$ (and also on $G$).     
Assume that $G$ is unramified at a rational $p$ and  the $p$-component $K_p$ of $K$ is a compact open subgroup of $G(\Z_p)$ which 
contains an Iwahori subgroup. 
There exists a compact open subgroup $K^p$ of $G(\A^{(p)}_f)$ such that $K=K^p\times K_p$.  
Note that $F_{f,p}$ is an Iwahori fixed vector of $\pi_{f,p}$.  
For $\alpha \in G(\Q)\cap T(\Q_p)/T(\Z_p)$, thus we obtain
\[
\begin{array}{rl}
T(\alpha)F_f(g)&=T(\alpha)F_f(g_\infty)\\
\\
&=T(\alpha)J^{-1}_\lambda(g_\infty,I)F(g_\infty  I)\\
\\
&=n(\alpha)_\lambda \ds\sum_iJ^{-1}_\lambda(g_\infty,I)
J^{-1}_\lambda(\alpha_{i,1},g_\infty I)F(\alpha_{i,1}g_\infty  I)\\
\\
&=n(\alpha)_\lambda \ds\sum_i
J^{-1}_\lambda(\alpha_{i,1}g_\infty, I)F(\alpha_{i,1}g_\infty  I)\\
\\
&=n(\alpha)_\lambda \ds\sum_i F_f(\alpha_{i,1}g_\infty).
\end{array}
\]
Note that 
$$G(\Q)Z_G(\R)^+\alpha_{i,1} g_\infty K=G(\Q)Z_G(\R)^+\alpha_{i} g_\infty K=
G(\Q)(g_\infty Z_G(\R)^+ \times \alpha^{-1}_i K).$$
Hence we have 
\[
\begin{array}{rl}
T(\alpha)F_f(g)&=n(\alpha)_\lambda \ds\sum_i F_f(g_\infty \alpha^{-1}_i)\\
\\
&=n(\alpha)_\lambda \ds\sum_{h\in K\alpha^{-1}K/K} F_f(g_\infty h)\\
\\
&=n(\alpha)_\lambda \ds\int_{K}F_f(g_\infty hg_f) dg_f\\
\\
&=n(\alpha)_\lambda \ds\int_{K\alpha^{-1}K}F_f(g_\infty g_f) dg_f\\
\\
&=n(\alpha)_\lambda \ds\int_{G(\A_f)}T(\alpha^{-1})_K F_f(g_\infty g_f) dg_f\\
\\
&=n(\alpha)_\lambda [K_p\alpha^{-1}K_p] F_{f,p}  
\end{array}
\]
where $dg_f$ is the Haar measure on $G(\A_f)$ so that vol$(K)=1$ and 
$T(\alpha^{-1})_K:=[K_p\alpha^{-1}K_p]\otimes 1_{K^p}$. 

We can naturally define the notion of positive or negative elements in the global setting. 
In the global setting,
we usually consider 
Hecke operators represented by  negative elements.
By the arguments above,
the computation of such global Hecke operators
is reduced to that of Hecke operators 
associated to positive elements in the local setting.     
In the next section, 
we will give examples of $p$-stabilized forms in 
the global setting. 

\section{$p$-stabilized forms}
\subsection{$GL_2$-case}In this subsection we will refer \cite{shimura} as a basic reference. 
Put $G=GL_2$. 
In this case,
we have $G^{{\rm der}}=SL_2,\ K=SO(2)(\R),\ K_\C\simeq \C^\times$ and the corresponding Hermitian symmetric space is the 
upper half-plane $\mathbb{H}=\{z\in \C|\ {\rm Im}(z)>0\}$. The automorphic factor is given by
$J(\gamma,z)=cz+d$ for 
$\gamma=
\left(\begin{array}{cc}
a& b\\
c& d
\end{array}
\right)
\in SL_2(\R)$ and  for an integer $k\ge 1$ we define the algebraic representation $\sigma_k:K_\C\lra \C^\times,z\mapsto z^k$.  
For an integer $N\ge 1$, we define the congruence subgroup $\Gamma_0(N)$ (respectively $\Gamma_1(N)$) to be 
the group consisting of the elements 
$g=\left(\begin{array}{cc}
a& b\\
c& d
\end{array}
\right)\in SL_2(\Z)$ such that $c\equiv 0 \ {\rm mod}\ N$ (respectively $a-1\equiv c\equiv 0 \ {\rm mod}\ N$). 
For an integer $N\ge 1$ and a Dirichlet character $\chi:\Gamma_0(N) \lra \C^\times$ so that $\chi |_{\Gamma_1(N)}=1$, 
we define $M_k(\Gamma_0(N),\chi):=M_{\sigma_k}(\Gamma_0(N),\chi)$. 
For a prime $p{\not |} N$, we define the action of Hecke operator $T_p=[\Gamma_0(N){\rm diag}(1,p)
\Gamma_0(N)]$ on $M_k(\Gamma_0(N),\chi)$ by 
\[
\begin{array}{rl}
T_pf(z):=&p^{\frac{k}{2}-1}\ds\sum_{\alpha\in \Gamma_0(N)\backslash\Gamma_0(N){\rm diag}(1,p)
\Gamma_0(N)}((p^{-\frac{1}{2}}\alpha)^{-1}f)(z)\\
\\
=&p^{\frac{k}{2}-1}\ds\sum_{\alpha\in \Gamma_0(N)\backslash\Gamma_0(N){\rm diag}(1,p)
\Gamma_0(N)}j(p^{-\frac{1}{2}}\alpha,z)^{-k}f(p^{-\frac{1}{2}}\alpha z).
\end{array}
\] 
Note that $p^{-\frac{1}{2}}\alpha\in SL_2(\R)$. 
If $p\not | N$, then we define the action of $U_p=[\Gamma_0(pN){\rm diag}(1,p)
\Gamma_0(pN)]$ on $M_k(\Gamma_0(pN),\chi)$ by 
$$U_pf(z):=p^{\frac{k}{2}-1}\sum_{\alpha\in \Gamma_0(pN)\backslash\Gamma_0(pN){\rm diag}(1,p)
\Gamma_0(pN)}((p^{-\frac{1}{2}}\alpha)^{-1}f)(z).$$
Let $f\in M_k(\Gamma_0(N),\chi)$ be a normalized cusp form which is an eigenform for all $T_p,\ {p\not|} N$ with 
the eigenvalue $a_p(f)$. Let $\alpha_p, \beta_p$ be 
the Satake parameters at $p$ so that 
$a_p(f)=\alpha_p+\beta_p$ and $\alpha_p\beta_p=\chi(p)p^{k-1}$. 
Let $\pi_{f}$ be the automorphic representation associated to $f$ and 
$\pi_{f,p}$ the local component at $p$.  
Since $\pi_{f,p}$  is a principal series representation, there exist 
characters $\chi_i:\Q^\times_p\lra \C^\times,\ i=1,2$ such that $\pi_{f,p}\simeq \pi(\chi_1,\chi_2)$. 
We may put $\chi_1(p^{-1})p^{\frac{k-1}{2}}=\alpha_p,\ 
\chi_2(p^{-1})p^{\frac{k-1}{2}}=\beta_p$.  
We define $p$-stabilized forms as follows:
$$f_{\alpha_p}(z):=f(z)-\beta_pf(pz),\ f_{\beta_p}(z):=f(z)-\alpha_pf(pz).$$
Then $f_{\alpha_p},\ f_{\beta_p}\in M_k(\Gamma_0(pN),\chi)$ and these are 
Hecke eigenforms for $T_\ell,\ \ell{\not|}pN$ and 
$U_p$-eigenforms with the eigenvalues $\alpha_p,\beta_p$ respectively. 

This can be checked via local computation as follows.  
By using the strong approximation theorem,
it is easy to see that 
$F_{f(pz)}(g)=F_f(g\cdot{\rm diag}(1,p))$. 
Hence the local component of $F_{f_{\alpha_p}}$ (respectively $F_{f_{\beta_p}}$) at $p$ 
corresponds to $\chi_1(p)f_1$ (respectively $\chi_2(p)f_2$) of Section 5. By the computation of 
Section 5 again, we have 
$$U_p f_{\alpha_p}(z)=p^{\frac{k}{2}-1}p^{\frac{1}{2}}\chi^{-1}_2(p) f_{\alpha_p}(z)=
\alpha_p f_{\alpha_p}(z)$$
and 
$$U_p f_{\beta_p}(z)=p^{\frac{k}{2}-1}p^{\frac{1}{2}}\chi^{-1}_1(p) f_{\beta_p}(z)=
\beta_p f_{\beta_p}(z).$$ 
We say the above $f$ is $p$-ordinary if ord$_p(a_p(f))=0$. 
If so is $f$, then we may assume that ord$_p(\alpha_p)=0$ and ord$_p(\beta_p)>0$.
The above computations show us that $f_{\alpha_p}$ can be embedded into a Hida family 
if $f$ is $p$-ordinary.  
Note that 
$f_{\beta_p}$ can be embedded into a Coleman family (cf. \cite{col}).
\subsection{$GSp_4$-case}
Let $\nu:G=GSp_4\lra GL_1$ be the similitude character. 
Put $Sp_4:={\rm Ker}\nu$. 
In this case, we have $G^{{\rm der}}=Sp_4$,\ 
$K=\Bigg\{\left(\begin{array}{cc}
A& B\\
-B& A
\end{array}
\right)\in {\rm Sp}_4(\R)\Bigg\}\simeq {\rm U}(2)(\R), K_\C\simeq {\rm GL}_2(\C)$ 
and the corresponding Hermitian symmetric space is the Siegel 
upper half-plane $\mathcal{H}_2=\{Z\in M_2(\C)|\ {}^tZ=Z,\  {\rm Im}(Z)>0\}$.  
For a pair of positive integers $\underline{k}=(k_1,k_2)$ such that $k_1\ge k_2$, we define the 
algebraic representation $\lambda_{\underline{k}}$ of ${\rm GL}_2(\C)$ by 
$$\lambda_{\underline{k}}={\rm Sym}^{k_1-k_2}{\rm St}_2\otimes {\rm det}^{k_2} {\rm St}_2$$ 
where ${\rm St}_2$ is the standard representation of dimension 2 over $\C$. 
Then the corresponding automorphic factor is defined by 
$$J_{\underline{k}}(\gamma,Z)=\lambda_{\underline{k}}(CZ+D)$$ 
 for 
$\gamma=
\left(\begin{array}{cc}
A& B\\
C& D
\end{array}
\right)
\in {\rm Sp}_4(\R)$ and $Z\in \mathcal{H}_2$. 
For an integer $N\ge 1$, we define a principal congruence subgroup $\Gamma(N)$ to be 
the group consisting of the elements $g\in {\rm Sp}_4(\Z)$ such that $g\equiv 1 \ {\rm mod}\ N$. 
We also define the level of $\Gamma$ to be the minimal $N$ satisfying (i)  
$g\Gamma g^{-1}$ contains $\Gamma(N)$ for some $g\in {\rm GSp}_4(\Q)$ and (ii)  all 
divisors of the denominator or the numerator of the entries of $g$ divide $N$. 
If $\Gamma$ is of level $N$, then 
the closure $\overline{\Gamma}$ in ${\rm Sp}_4(\A_f)$ satisfies that 
the $p$-component of $\overline{\Gamma}$ is ${\rm Sp}_4(\Z_p)$ for all $p{\not|}N$. 
For a parabolic subgroup $P$, let $I_P=\Gamma_P(p)$ be the group consisting of the elements $g\in {\rm Sp}_4(\Z)$ such that $(g\ {\rm mod}\ p)\in 
P(\F_p)$. For a discrete subgroup $\Gamma$ of level $N$, put $\Gamma_{I_P}:=\Gamma\cap I_P$. 

For a discrete subgroup $\Gamma$ of level $N$ and 
 a Dirichlet character $\chi:\Gamma \lra \C^\times$, 
we define $M_{\underline{k}}(\Gamma,\chi):=M_{\lambda_{\underline{k}}}(\Gamma,\chi)$. 
A function of this space is called a Siegel modular form of weight $\underline{k}$ and 
level $\Gamma$ with a character $\chi$. 
If $k=k_1=k_2$, put $M_{k}(\Gamma,\chi):=M_{(k,k)}(\Gamma,\chi)$ for short. 
For a prime $p\not | N$ and $t'_1={\rm diag}(1,1,p,p),\ t'_2={\rm diag}(1,p,p,p^2)$, we define the action of two Hecke operators 
$T_{p,i}=[\Gamma t'_1\Gamma],\ i=1,2$  on $M_k:=\ds\bigcup_{(\Gamma,\chi)} M_k(\Gamma,\chi)$ by 
\[
\begin{array}{rl}
T_{p,i}f(Z):=&p^{i(\frac{k_1+k_2}{2}-3)}\ds\sum_{\alpha\in \Gamma\backslash\Gamma t'_i
\Gamma}((\nu(t'_i)^{-\frac{1}{2}}\alpha)^{-1}f)(Z)  \\
\\
=&p^{i(\frac{k_1+k_2}{2}-3)}\ds\sum_{\alpha\in \Gamma\backslash\Gamma t'_i
\Gamma}J_{\underline{k}}(\nu(t'_i)^{-\frac{1}{2}}\alpha,Z)^{-1}f(\nu(t'_i)^{-\frac{1}{2}}\alpha Z).
\end{array}
\]
Then $T_{p,i}f\in M_k$, but in general, these operators do not preserve $M_{\underline{k}}(\Gamma,\chi)$ (see Lemma 3.1 of \cite{mt}).
Note that $\nu(t'_i)^{-\frac{1}{2}}\alpha\in  {\rm Sp}_4(\R)$. 
If $p{\not |} N$, then we define the action of $U^{P,{\rm global}}_{p,i}=[\Gamma_{I_P} t'_i\Gamma_{I_P}],\ i=1,2$ on $M_k(\Gamma_{I_P},\chi)$ by 
\begin{equation}\label{UpforGSp4}
U^{P,{\rm global}}_{p,i}f(Z):=p^{i(\frac{k_1+k_2}{2}-3)}\sum_{\alpha\in \Gamma_{I,P}\backslash\Gamma_{I_P} t'_i\Gamma_{I_P}}((\nu(t'_i)^{-\frac{1}{2}}\alpha)^{-1}f)(Z).
\end{equation}

In what follows, we will discuss about a $p$-stabilized form of Saito-Kurokawa lift. 
Let $\Gamma_0(N)$ be the subgroup of ${\rm Sp}_4(\Z)$ consisting of the elements 
$g=
\left(\begin{array}{cc}
A_g& B_g\\
C_g& D_g
\end{array}\right)$ 
such that $C_g\equiv 0_2$ mod $N$. Fix a finite character $\chi:(\Z/N\Z)^\times\longrightarrow \mathbb{C}^\times$. 
Then we define (and denote it by $\chi$ again)
the character on $\Gamma_0(N)$ associated to $\chi$ by 
$\chi(g):=\chi(\det(D_g))$ for $g\in \Gamma_0(N)$.
For a normalized elliptic newform $f$ of weight $2k-2\ge 2$ and level $N$ with the character $\chi^2$, 
there exists a cusp form $F=SK(f)$ in $S_{k-1}(\Gamma_0(N),\chi)$ by \cite{ibu} so that 
$F$ is an eigenform for all $T_{p,i},\ p{\not |} N$ with 
the eigenvalue $a_{p,i}(F)$ and these eigenvalues are written as  
$$a_{p,1}(F)=\chi(p)(p^{k-1}+p^{k-2})+a_p(f),$$
$$a_{p,2}(F)=a_p(f)\chi(p)(p^{k-2}+p^{k-3})+2\chi^2(p)p^{2k-4}-(p^2+1)\chi^2(p)p^{2k-6}$$
where $a_p(f)$ is the eigenvalue of $f$ for $T_p$ in Section 9.1. 
Let $\pi_{f}$ (respectively $\Pi_F$) be the automorphic representation associated to $f$ (respectively $F=SK(f)$) and 
$\pi_{f,p}$ (respectively $\Pi_{F,p}$) the local component at $p$. Let $\chi':\A^\times \lra \C^\times$ be 
the character corresponding to $\chi$ and denote by $\chi'_p$ its local component at $p$. 
Since $\pi_{f,p}$  is a principal series representation, there exist 
characters $\chi_i:\Q^\times_p\lra \C^\times,\ i=1,2$ such that $\pi_{f,p}\simeq \pi(\chi_1,\chi_2)$ 
with the central character $\chi_1\chi_2=\chi'^{-2}_p$. 
Hence we have $\chi_1\chi_2(p)=\chi'^2_p(p^{-1})$.  
Then by \cite{rs}, we have $\Pi_{F,p}\simeq \chi_1\chi' 1_{{\rm GL}_2}\rtimes \chi^{-1}_1\chi'^{-2}_p$ which is  
the case IIb in Section 7.1.2. Then the eigenvalues of $U^B_{p,i},i=1,2$ in the local setting are 
$$(p^{\frac{3}{2}}\delta,p^{2}\delta\beta),(p^{\frac{3}{2}}\alpha,
p^{2}\alpha\beta),(p^{\frac{3}{2}}\beta,p^{2}\beta\delta),
(p^{\frac{3}{2}}\beta,p^{2}\beta\alpha)$$
where $\alpha=\chi_1(p^{-1}),\ \beta=p^{\frac{1}{2}}\chi'(p),\ \gamma
=p^{-\frac{1}{2}}\chi'(p)$ and $\delta=\chi_2(p^{-1})$. 
We may put $\chi_1(p^{-1})p^{\frac{2k-3}{2}}=\alpha_p(f),\ 
\chi_2(p^{-1})p^{\frac{2k-3}{2}}=\beta_p(f)$ where $\alpha_p(f),\beta_p(f)$ are eigenvalues of 
$T_p$ for $f$. 
From (\ref{UpforGSp4}) and the observation in Section 7 (7.1-(3) and 7.1.2 Case IIb), 
\begin{equation}
(\beta_p-U^{B,{\rm global}}_{p,1})(\chi'(p)p^{k-1}-U^{B,{\rm global}}_{p,1})F
\end{equation}
is a $p$-stabilized form with respect to $\widehat{B}=B$ with the eigenvalues $(\alpha_p,p^{k-2}\alpha_p\chi'(p))$ for 
$U^{B,{\rm global}}_{p,i},i=1,2$. 
If we take the  Fourier expansion $F=\ds\sum_{T\in {\rm Sym}^2(\Z)_{>0}}a(T)e^{2\pi\sqrt{-1}{\rm tr}(TZ)}$ where 
${\rm Sym}^2(\Z)_{>0}$ is the subset of $M_2(\Q)$ consisting of all symmetric matrices which are positive and semi-integral, 
then we have  
$$(\beta_p-U^{B,{\rm global}}_{p,1})(\chi'(p)p^{k-1}-U^{B,{\rm global}}_{p,1})F=$$
\begin{equation}
\ds\sum_{T>0}\Big(p^{k-1}\chi'(p)\beta_p a(T)-(\chi'(p)p^{k-1}+\beta_p) a(pT)+ a(p^2T)\Big)e^{2\pi\sqrt{-1}{\rm tr}(TZ)},
\end{equation}
since we can choose a complete system of representatives for $\Gamma_{I,B}\backslash\Gamma_{I_B} t'_1\Gamma_{I_B}$ 
to be the same as  in the case Siegel parabolic and in that case we can compute the action easily (cf. \cite{bocherer}).  

Similarly, 
$(\beta_p-U^{P,{\rm global}}_{p,1})(\chi'(p)p^{k-1}-U^{P,{\rm global}}_{p,1})F$
is also a $p$-stabilized form with respect to $\widehat{P}=Q$ with the eigenvalues $\alpha_p$ for 
$U^{P,{\rm global}}_{p,1}$ (see 7.2-(3) and 7.2.2 Case IIb).  
As explained above, we have the same expansion. 
Therefore  we have 
$$(\beta_p-U^{B,{\rm global}}_{p,1})(\chi'(p)p^{k-1}-U^{B,{\rm global}}_{p,1})F=(\beta_p-U^{P,{\rm global}}_{p,1})(\chi'(p)p^{k-1}-U^{P,{\rm global}}_{p,1})F$$ 
for $t'_1$. 
If $k\ge 2$, $F$ has $p$-integral (algebraic) coefficients, then  so does 
$(\beta_p-U^{B,{\rm global}}_{p,1})(\chi'(p)p^{k-1}-U^{B,{\rm global}}_{p,1})F$. 
We further assume that ${\rm ord}_p(\alpha_p)=0$ (this implies  ${\rm ord}_p(\beta_p)>0$). 
Then by using the relation of Hecke operators (cf. p.228, Example 4.2.10 of \cite{And}),
we see that ${\rm ord}_p(a(p^2T_0))=0$ if ${\rm ord}_p(a(T_0))=0$ for some $T_0\in {\rm Sym}^2(\Z)_{>0}$. 
Hence under this assumption (the existence of such $T_0$), 
$(\beta_p-U^{B,{\rm global}}_{p,1})(\chi'(p)p^{k-1}-U^{B,{\rm global}}_{p,1})F$ is not only preserved $p$-integrality, but also  
non-vanishing modulo $p$.



\begin{thebibliography}{99}
{
\bibitem{And}A. N. Andrianov, Quadratic forms and Hecke operators, Springer Berlin 1987. 
\bibitem{bocherer}S. B\"ocherer, Siegfried On the Hecke operator $U(p)$ . 
With an appendix by Ralf Schmidt. J. Math. Kyoto Univ.  45  (2005),  no. 4, 807-829. 
\bibitem{Borel}A. Borel, Automorphic L -functions. Automorphic forms, representations and L -functions (Proc. Sympos. Pure Math., Oregon State Univ., Corvallis, Ore., 1977), Part 2, pp. 27-61, Proc. Sympos. Pure Math., XXXIII, Amer. Math. Soc., Providence, R.I., 1979.
\bibitem{BJ}A. Borel and H. Jacquet, 
 Automorphic forms and automorphic representations. 
With a supplement "On the notion of an automorphic representation'' by R. P. Langlands. Proc. Sympos. Pure Math., XXXIII, Automorphic forms, representations and L -functions (Proc. Sympos. Pure Math., Oregon State Univ., Corvallis, Ore., 1977), Part 1, pp. 189-207. 
\bibitem{bk}C. J. Bushnell and P. C. Kutzko, Smooth representations of reductive p -adic groups: structure theory via types. Proc. London Math. Soc. (3) 77 (1998), no. 3, 582-634.
\bibitem{cas-note}W. Casselman, Introduction to admissible representations of $p$-adic groups, available at 
his homepage.  
\bibitem{col}Robert F. Coleman, Classical and overconvergent modular forms, Invent. Math. 124 (1996), no. 1-3, 215-241.
\bibitem{garrett}P. Garrett, Representations with Iwahori-fixed vectors, note,  available at 
http://www.math.umn.edu/~garrett/m/v/.
\bibitem{Goresky}M. Goresky, Compactifications and cohomology of modular varieties. Harmonic analysis, the trace formula, and Shimura varieties, 551-582, Clay Math. Proc., 4, Amer. Math. Soc., Providence, RI, 2005.
\bibitem{Harish}
R.~Howe.
\newblock {\em Harish-{C}handra homomorphisms for {${p}$}-adic groups},
  volume~59 of {\em CBMS Regional Conference Series in Mathematics}.
\newblock Published for the Conference Board of the Mathematical Sciences,
  Washington, DC, 1985.
\newblock With the collaboration of A. Moy.
\bibitem{ibu}T. Ibukiyama, Saito-Kurokawa liftings of level N and practical construction of Jacobi forms, 
Kyoto J. Math. Volume 52, Number 1 (2012), 141-178.
\bibitem{jacquet}H. Jacquet, Sur les representations des groupes reductifs $p$-adiques. 
C. R. Acad. Sci. Paris Ser. A-B 280 (1975), Aii, A1271-A1272. 
\bibitem{keys}D. Keys, Principal series representations of special unitary groups over local fields. Compositio Math. 51 (1984), no. 1, 115-130.
\bibitem{Mazur}B.~Mazur,
\newblock An ``infinite fern'' in the universal deformation space of {G}alois
  representations.
\newblock {\em Collect. Math.}, 48(1-2):155--193, 1997.
\newblock Journ{\'e}es Arithm{\'e}tiques (Barcelona, 1995).
\bibitem{Milne}J.~S. Milne,
\newblock Introduction to {S}himura varieties.
\newblock In {\em Harmonic analysis, the trace formula, and {S}himura
  varieties}, volume~4 of {\em Clay Math. Proc.}, pages 265--378. Amer. Math.
  Soc., Providence, RI, 2005.
\bibitem{rs}B. Roberts and R. Schmidt, 
Local newforms for $GSp(4)$. Lecture Notes in Mathematics, 1918. Springer, Berlin, 2007. viii+307 pp.
\bibitem{mt}R. Salvati Manni and J. Top, Cusp forms of weight 2 for the group $\Gamma(4,8)$, Amer. J. Math. 115 (1993), 455-486. 
\bibitem{Satake}I. Satake, Algebraic structures of symmetric domains. Kano Memorial Lectures, 4. Iwanami Shoten, Tokyo; Princeton University Press, Princeton, N.J., 1980. xvi+321 pp.
\bibitem{shimura}G. Shimura, Introduction to the arithmetic theory of automorphic functions. Reprint of the 1971 original. Publications of the Mathematical Society of Japan, 11. Kano Memorial Lectures, 1. Princeton University Press, Princeton, NJ, 1994.
\bibitem{skinner&urban}C. Skinner and E. Urban,  Sur les d\'eformations $p$-adiques de certaines representations automorphes. J. Inst. Math. Jussieu 5 (2006), no. 4, 629-698.
\bibitem{til}J. Tilouine, Nearly ordinary rank four Galois representations and $p$-adic Siegel modular forms. With an appendix by Don Blasius. Compos. Math. 142 (2006), no. 5, 1122-1156. 
}
\end{thebibliography}
\end{document}